\documentclass{amsart}
\usepackage{amssymb, amsmath, latexsym, amscd}
\usepackage{epsfig}

\newcommand{\lbl}{\label}
\newcommand{\biblbl}{\bibitem}
\newtheorem{theorem}{Theorem}
\newtheorem{definition}[theorem]{Definition}
\newtheorem{corollary}[theorem]{Corollary}
\newtheorem{lemma}[theorem]{Lemma}
\newtheorem{proposition}[theorem]{Proposition}

\newtheorem{remark}{Remark} 
\newtheorem{remarks}{Remarks} 
\newtheorem{thm}{Theorem}

\newcommand{\Z}{\mathbb{Z}}
\newcommand{\N}{\mathbb{N}}
\newcommand{\Q}{\mathbb{Q}}
\newcommand{\R}{\mathbb{R}}
\newcommand{\rr}{\mathcal{R}}

\newcommand{\F}{\mathcal{F}}
\newcommand{\h}{\mathcal{H}}
\newcommand{\K}{\mathcal{K}}
\newcommand{\T}{\mathcal{T}}

\newcommand{\I}{\mathcal {I}}

\newcommand{\C}{\mathcal{C}}

\newcommand{\bo}{\partial}
\newcommand{\G}{\Gamma}
\newcommand{\down}{\downarrow}
\newcommand{\up}{\uparrow}

\begin{document}
\title{Grope cobordism of classical knots}

\author[J. Conant]{James Conant}
\address{Department of Mathematics\\
         Cornell University\\
         Ithaca, NY 14853-4201}
\email{jconant@math.cornell.edu}

\author[P. Teichner]{Peter Teichner}
\address{Department of Mathematics,
         University of California in San Diego,
         9500 Gilman Dr, La Jolla, CA 92093-0112}
\email{teichner@euclid.ucsd.edu}

\thanks{The first author is partially supported by
NSF VIGRE grant DMS-9983660, and the second by NSF grant DMS0072775.} 

\begin{abstract}
{Motivated by the lower central series of a group, we define
 the notion of a {\em grope cobordism} between two knots in a
$3$-manifold. Just like an iterated group commutator, each grope cobordism has a type that
can be described by a rooted unitrivalent tree. By filtering these trees in different
ways, we show how the Goussarov-Habiro approach to finite type invariants of knots is
closely related to our notion of grope cobordism. Thus our results can be
viewed as a geometric interpretation of finite type invariants.

The derived commutator series of a group also has a 3-dimensional analogy, namely knots
modulo {\em symmetric} grope cobordism. On one hand this theory maps onto the
usual Vassiliev theory and on the other hand it maps onto the Cochran-Orr-Teichner
filtration of the knot concordance group, via symmetric grope cobordism in 4-space. In
particular, the graded theory contains information on finite type invariants
(with degree~$h$ terms mapping to Vassiliev degree $2^h$), Blanchfield forms
or S-equivalence at $h=2$, Casson-Gordon invariants at $h=3$, and for $h=4$
one finds the new von Neumann signatures of a knot. }
\end{abstract}

\maketitle

\section{Introduction}

A modern perspective on $3$-manifolds is through topological quantum field
theory, following ideas of Jones, Witten and many others. These have inspired
tremendous activity but so far have not contributed much to the
topological understanding of $3$-manifolds. In particular, the
Vassiliev-Goussarov theory of finite type invariants of knots, which in
some sense gives universal quantum knot invariants, has developed a
fascinating life quite independent of the rest of geometric topology.
Which low-dimensional topologist hasn't been inspired by the appearance
of unitrivalent graphs in the enumeration of these finite type
invariants? These graphs ultimately derive from the Feynman rules
associated to perturbative Chern-Simons theory, and the residue of Gauge
symmetry introduces certain relations on the diagrams, now known as {\em
antisymmetry-} and {\em Jacobi-} (or IHX-) relations. On the other hand,
it is well known that rooted unitrivalent {\em trees} can be used to
label iterated (non-associative) operations, and that the above relations
arise exactly for Lie algebras. In our context the most interesting Lie
algebras arise from a group $G$ by first considering its {\em lower
central series} $G_c$ defined inductively by the iterated commutators
$$ 
G_2:=[G,G] \text{ and } G_c:=[G,G_{c-1}] \text{ for } c>2.
$$ 
Then $L(G):=\bigoplus_c G_c$ is a Lie algebra with group multiplication as addition and group
commutators as Lie bracket.

In this paper we shall give a geometric implementation of iterated
commutators in fundamental groups $G$ via the notion of a {\em grope cobordism} between two
knots in a $3$-manifold. From the above point of view, one should rather think of the
associated graded Lie algebra $L(G)$ and thus it is not surprising that the notion of grope
cobordism is closely related to finite type knot invariants. We shall make this statement
 precise and hence our results can be viewed as the long desired {\em
geometric} interpretation of finite type knot invariants. 

This relation between grope cobordism and finite type invariants was first
announced by Habiro at the very end of his landmark paper
\cite{h2}. Without providing proofs, he correctly announces a version of
Theorem~\ref{capped} below, but makes an incorrect assertion about
(uncapped) grope cobordism. The correct statement is our main result,
Theorem~\ref{main}. Our proofs of these theorems rely heavily on Habiro's
work.

Other geometric interpretations of finite type invariants include Stanford's 
beautiful work \cite{s} on the relationship with the lower central
series of pure braid groups $PB_n$. Stanford shows that two knots in 3-space have the same
finite type invariants of degree~$<c$ if and only if they differ by a finite sequence of
operations as follows: Grab any number $n$ of strands of one knot and tie them into a pure
braid in the $c$-th term of the lower central series of $PB_n$.

The first relation between finite type invariants and gropes was
announced by Kalfagianni and Lin in \cite{kl}. Their notion is very different from ours since
they consider gropes in 3-space whose first stage bounds a knot and is embedded with free
complementary fundamental group. However, arbitrary intersections are allowed among the higher
grope stages. In that context the precise relation between grope class and Vassiliev degree
is not understood and only a logarithmic estimate is given in \cite{kl}.
In his thesis \cite{c}, the first author discovered a more precise relationship between
finite-type invariants  and gropes. There he proved that a knot bounding an
{\em embedded} grope of class $c$ in 3-space must have vanishing finite type invariants up to
$\lceil c/2 \rceil$, and that this bound is the best possible. The methods of the thesis are applied in the short note \cite{newconant} to get a similar result for gropes with more than one boundary component. This result is an 
ingredient of the proof of Theorem~\ref{main} below.

\subsection{A geometric interpretation of group commutators}
We first want to motivate the use of gropes from scratch, without any
reference to quantum invariants. Recall that the fundamental group of an
arbitrary topological space $X$ consists of continuous maps of the circle
$S^1$ into $X$, modulo homotopy (i.e.\
$1$-parameter families of continuous maps). Quite analogously, classical
knot theory studies smooth {\em embeddings} of a circle into 
3-space, modulo isotopy (i.e.\ $1$-parameter families of embeddings). 

\begin{figure}[ht]
\begin{center}
\epsfig{file=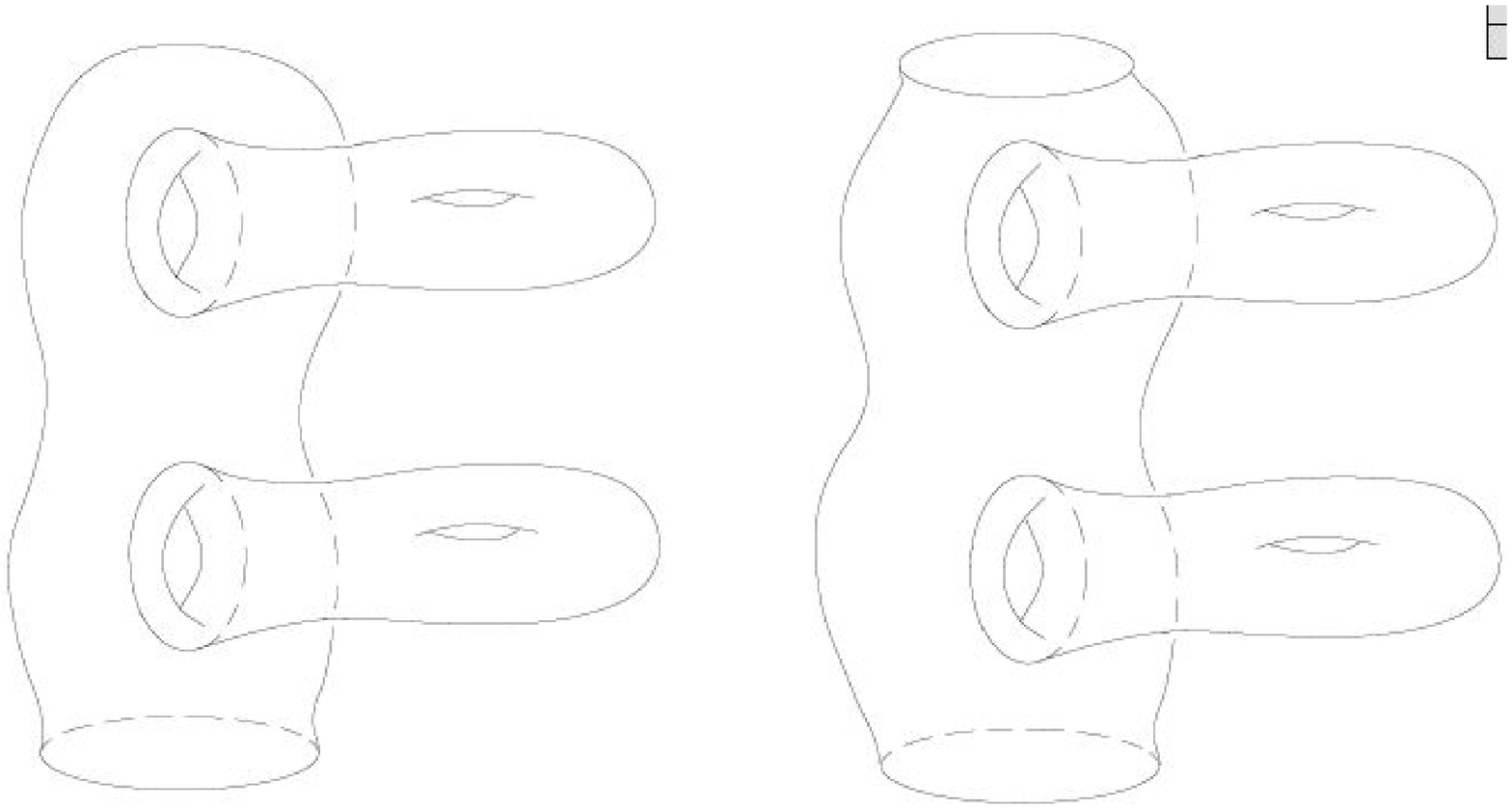,%
            height = 4cm}
\end{center}
\caption{Gropes of class 3, with one and two boundary circles.}
\lbl{fig:grope}
\end{figure}

Recall further that a continuous map $S^1\to X$
represents the trivial element in the fundamental group $\pi_1X$ if and
only if it extends to a map of the disk. Moreover,
$\phi$ represents a commutator in
$\pi_1X$ if and only if it extends to a map of a surface (i.e.\ of
a compact oriented 2-manifold with boundary $S^1$). The first statement has
a straightforward analogy in knot theory: A knot is trivial if and
only if it extends to an embedding of the disk into 3-space. However, every
knot ``is a commutator'' in the sense that it bounds a {\em Seifert
surface}, i.e.\ an embedded surface in 3-space. 
Thus all of knot theory is created by the difference between a surface and a disk.
The new idea is to filter this difference by introducing a concept into
knot theory which is the embedded analogue of iterated commutators in group theory.
Namely, there are certain finite 2-complexes (built out of iterated surface
stages) called {\em gropes} by Cannon
\cite{ca}, with the following defining property: $S^1\to X$ represents
an element in the $c$-th term of the lower central series
of $\pi_1X$ if and only if it extends to a continuous map of a {\em grope
of class~$c$}. By construction, such gropes have a single circle as their boundary, but one
can also consider gropes with more boundary circles as in Figure~\ref{fig:grope}.

Gropes, therefore, are not quite manifolds but the singularities that arise
are of a very simple type, so that these 2-complexes are in some sense the
next easiest thing after surfaces. Two sentences on the history of the use
of gropes in mathematics are in order, compare \cite[Sec.2.11]{fq}. Their
inventor Stan'ko worked in high-dimensional topology, and so did Edwards
and Cannon who developed gropes further. Bob Edwards suggested their
relevance for topological 4-manifolds, where they were used extensively,
see \cite{fq}, \cite{ft1}, or \cite{ft2}. It is this application that
seems to have created a certain angst about studying gropes, so we should
point out that the only really difficult part in 4 dimensions is the use
of {\em infinite constructions}, i.e.\ when the class of the grope goes to
infinity. 

One of the purposes of this paper is to show how simple and useful (finite) gropes are when
embedded into 3-space.

\subsection{Grope cobordism of knots in 3-space}
The idea behind a {\em grope cobordism} is to 
filter the difference between a surface and a disk
in 3-space. The following definition should be thought of as a 3-dimensional {\em
embedded} analogue of the lower central series of a group. Let $\K$ be the set of oriented
knot types, i.e. isotopy classes of oriented knots in 3-space.

\begin{definition} \lbl{def:grope cobordism}
Two knot types $K_1,K_2\in\K$ are
{\em grope cobordant} of class~$c$, if there is an embedded grope of class~$c$ 
(the {\em grope cobordism}) in 3-space such that its two boundary components represent $K_1$
and $K_2$.
\end{definition}

\begin{figure}[ht]
\begin{center}
\epsfig{file=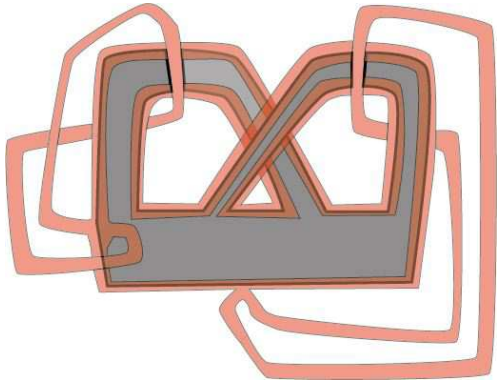,height=5cm}
\caption{A class three grope cobordism.}\label{fig:realisticgrope}
\end{center}
\end{figure}

At first glance, gropes don't appear to embed in an interesting way in
$3$-space. However,
 since every grope cobordism has a 1-dimensional
spine, it can then be isotoped into the neighborhood of a
1-complex. As a consequence, grope cobordisms abound in 3-space!  An example of such a grope cobordism of
class three is given in Figure~\ref{fig:realisticgrope}. This is an embedded version of the grope on
the right of Figure~\ref{fig:grope}, except that all surface stages are of genus one. 
The genus one surface with two boundary components is the thin, partially transparent surface.
One symplectic basis element, the core of one of the thin bands, is glued to the boundary
of the thicker genus one surface.
It is
important to point out that the two boundary components of a grope cobordism may link in an
arbitrary way, but that we do not record this information.

It turns out (Lemma~\ref{lem:equivalence relation}) that the relation of grope cobordism is
in fact an equivalence relation (for each fixed class ~$c$) on the set $\K$ of knot types.
This is why we were careful to talk about {\em knot types} rather than actual knots.
Moreover, the resulting quotients are extremely interesting abelian groups under connected
sum. Before explaining these groups in detail, we want to point out a way to directly relate
to finite type knot invariants a l\'a Vassiliev, see
\cite{v} or \cite{bn}. For that purpose, we have to consider {\em capped} gropes which are
gropes with disks (the {\em caps}) as their top surface stages. 

If two knots cobound an {\em embedded capped} grope then they are isotopic because
the caps can be used to surger the grope cobordism into an annulus. Thus in order to
get an interesting notion of {\em capped grope cobordism}, we allow the
(disjointly embedded) caps to have intersections with the bottom stage of the embedded
grope.

\begin{theorem}\lbl{capped} Two oriented knots are capped grope cobordant
of class $c$ if and only if they share the same finite type invariants of
Vassiliev degree $<c$.
\end{theorem}

The proof of this result has two ingredients. One is Habiro's beautiful
translation of finite type invariants into his theory of tree claspers
\cite{h2}, and the other is our translation from tree claspers to capped
gropes given in Theorem~\ref{claspers-gropes}.

\begin{remarks} In Section~\ref{sec:grope refinement} we will
prove that it is sufficient to consider gropes which have genus
one in all stages except at the bottom. The genus at the bottom
is responsible for the transitivity of the grope cobordism
relation.

Another simplification is predicted by group theory: Since the lower central
series of a group is generated by commutators which are ``right-normed'', the
question arises as to whether (capped) grope cobordism is generated by the
corresponding {\em half-gropes}, see Figure~\ref{fig:half}. This question will 
be answered in the affirmative in Section~\ref{sec:IHX}.
\end{remarks}

Even though capped grope cobordism is very useful because of Theorem~\ref{capped}, the
analogy with group theory is more natural in the absence of
caps. Thus the question arises whether grope cobordism (without caps) can
 also be translated into the finite type theory. In order to explain how
this can be done, we first have to review an approach to finite type knot
invariants developed by Goussarov, Habiro and others.

\subsection{Finite type filtrations and grope cobordism}
Again the starting point are certain Feynman diagrams, i.e. unitrivalent graphs.
The main idea is to think of such graphs as operating on the space of knots as follows. 
Consider a unitrivalent graph $\G$ embedded in 3-space, with exactly its univalent
vertices on a knot $K$, its edges framed and each trivalent vertices cyclically
ordered. There is a procedure to replace
$\G$ by a framed link in the complement of $K$, with a copy of the Hopf link
at each edge of $\G$ and a copy of the Borromean rings at each trivalent
vertex. See Section~\ref{sec:claspersdefs}.
 Surgery on that link 
replaces the knot $K$ by a new knot type $K_\G$, the {\em surgery of $K$ along
$\G$}. 
In the simplest case where
$\G_0$ has a single edge, one recovers the original idea of a crossing change on a knot:
Surgery on a single Hopf link leads to the knot
$K_{\G_0}$ which differs from $K$ by a single crossing change. The next simplest case is
shown in Figure~\ref{fig:intro8}.

\begin{figure}[ht]
\begin{center}
\epsfig{file=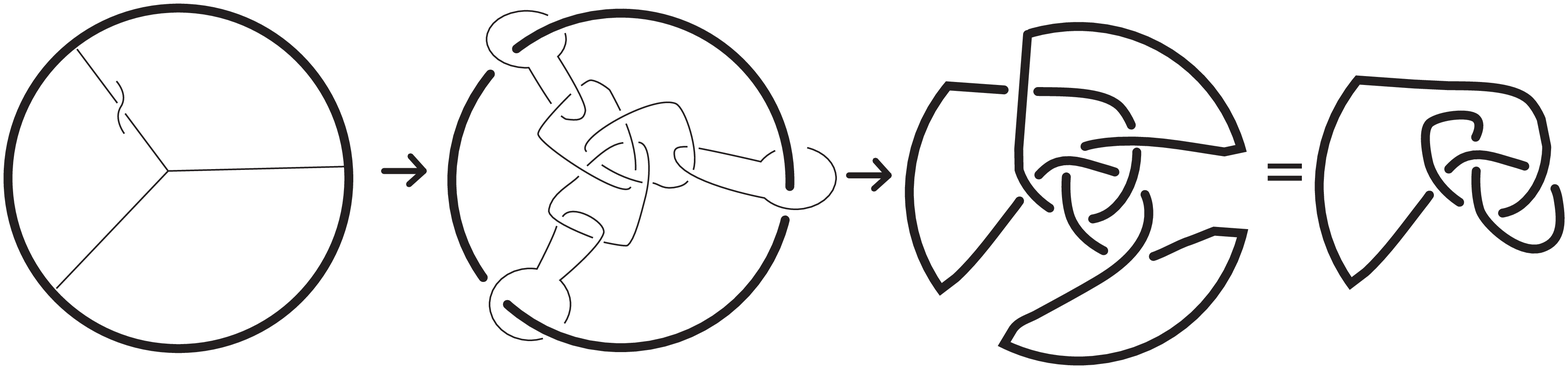,%
            height=2.6cm}
\end{center}
\caption{Surgering the unknot to the Figure~8 knot. The framing is the blackboard framing, except
at the indicated half-twist.}
\lbl{fig:intro8}
\end{figure}

Varying the embeddings and framings of a given graph, one obtains an infinite
class of operators on the set $\K$ of oriented knot types in 3-space, indexed by abstract
unitrivalent graphs.  Assume that each such graph $\Gamma$ is equipped with a {\em degree}
$\deg(\G)\in \N$. Then one obtains a descending filtration
$$
\K = \F^{\deg}_0 \supseteq \F^{\deg}_1 \supseteq \F^{\deg}_2 \supseteq \dots
$$ defined as follows: $\F^{\deg}_k$ consists of all knots that can be
obtained from the unknot by a finite sequence of surgeries along 
unitrivalent graphs $\G$ of degree $\deg(\G) \geq k$. There is also a natural
notion of the quotients $\K/\F_k^{\deg}$: These are defined to be the
equivalence classes of the equivalence relation on
$\K$ generated by surgeries along rooted unitrivalent graphs of degree
$\geq k$.

As an example, one can use the {\em Vassiliev degree}
$$ 
v(\G):=(\text{number of vertices of } \G)/2
$$
to obtain exactly the well known {\em Vassiliev filtration} used in
Theorem~\ref{capped}: The main theorem of
\cite{h2} states that two knots represent the same element in $\K/\F_k^v$ if and only if
they share the same Vassiliev invariants of degree $<k$. 
This follows from the fact that a surgery on a unitrivalent graph is the same as a {\em
simple clasper surgery}, compare Section~\ref{sec:claspers}. In this language,
Theorem~\ref{capped} can be reformulated as follows:

\begin{thm}\lbl{capped'} 
Two oriented knots in 3-space are capped grope cobordant
of class $c$ if and only if they represent the same element in $\K/\F_c^v$.
\end{thm}

The Vassiliev degree is also used as the degree in defining a version of
graph cohomology. Then it turns out that the differential in this chain
complex preserves another degree, namely the {\em loop degree}
$\ell(\G):=b_1(\G)$, the ``number'' of loops in
$\G$. Regardless of its relation to graph cohomology, one can use the loop
degree to obtain a second filtration $\F_k^{\ell}$ of the set of knot types.
It turns out that in our context the {\em grope degree}
$$ g(\G):= v(\G) +\ell(\G)
$$ 
is most relevant, leading to the precise uncapped analogue of Theorem~\ref{capped'}:

\begin{theorem}\lbl{main} Two oriented knots in 3-space are grope cobordant
of class $c$ if and only if they represent the same element in $\K/\F_c^g$.
\end{theorem}

\noindent {\bf Remarks:}
\begin{itemize}
\item The grope degree arises naturally as follows. Given a unitrivalent
graph
$\Gamma$, there exists a set of $\ell(\Gamma)$ edges such that cutting
these edges yields a trivalent tree. The grope degree of $\Gamma$ is
precisely the Vassiliev degree of this tree (which is the same as the
degree of the tree in the sense of the lower central series).  
\item The groups $\K/\F_c^v\otimes\Q$ are well known to be isomorphic to the
corresponding diagram spaces via the Kontsevich integral. In particular, the Kontsevich
integral is an invariant of capped grope cobordism. 
\item Garoufalidis and Rozansky
\cite{gr} have proven the remarkable result that the Kontsevich integral also preserves the
``loop filtration'' $\F^\ell_*$.
In particular, the Kontsevich integral also preserves the grope
filtration $\F^g_*$ (in the sense that it sends the $c$-th term $\F_c^g$ to a linear
combination of diagrams with grope degree $\geq c$). Hence the Kontsevich integral gives
obstructions to the existence of grope cobordisms. 
\item In fact, the groups $\K/\F_c^g\otimes\Q$ are isomorphic to the
corresponding diagram spaces (via the Kontsevich integral) just like for
the Vassiliev degree. This result will be explained in
\cite{ct}. It shows how interesting, yet understandable, the relation of grope cobordism in
3-space is. Moreover, it also gives a geometric interpretation of the Kontsevich integral!
\item Using the methods of \cite{h2} one shows that all quotients,
$\K/\F_c^v$ and $\K/\F_c^g$, are finitely generated abelian groups under
connected sum. Hence the same is true for the quotients of knots modulo
(capped) grope cobordism.
\item In the preceding theorems we are dividing out by
graphs which have grope degree {\em larger than} or equal to $c$. In
fact the theorems are also true if we only divide out by those of
exactly degree $c$. For the Vassiliev degree this is contained in \cite{h2}, for the grope
degree we shall give a proof in \cite{ct}.
\end{itemize}

 Theorem~\ref{main} will be proven by explaining the precise
relation between an embedded grope and the link obtained from a rooted 
unitrivalent graph. At the heart of the issue lies a well
known relation between Borromean rings and surfaces and more generally
between iterated Bing doublings of the Hopf link and gropes of higher class.
This relation has been used extensively in 4-dimensional topology and it has
also occurred previously in the study of Milnor invariants of links, see for
example \cite{co}.

\subsection{Gropes and claspers}
The following result is our main contribution to Theorems~\ref{capped} and
\ref{main}. It uses an obvious generalization of grope cobordism in 3-space to arbitrary
$3$-manifolds and also the language of {\em claspers} which makes the ``surgery on
unitrivalent graphs'' from the previous section more precise.

\begin{theorem}\lbl{claspers-gropes} Let $\T$ be a rooted trivalent tree
and let $M$ be a $3$-manifold.
\begin{itemize}
\item[(a)] Two knots are $\T$-grope cobordant in $M$ if and only if they are
related by a finite sequence of $\T$-clasper surgeries.
\item[(b)] Two knots are capped $\T$-grope cobordant in $M$ if and only if
they are related by a finite sequence of capped $\T$-clasper surgeries.
\end{itemize}
\end{theorem}

All the relevant definitions will be introduced in the next
sections. In particular, we shall explain the correspondence between rooted trees,
gropes and claspers.

\subsection{4-dimensional aspects}
By a result of Ng \cite{ng}, no finite type knot invariant but the Arf invariant is a
concordance invariant. The analogous result for links can be very well expressed in
terms of the loop degree and its relation to gropes. Even though our (grope) proof is new,
the following result seems well known to experts. Compare in particular the rational
analogue of Habegger and Masbaum \cite{hm} and the homological version of Levine \cite{l}.

\begin{theorem} \lbl{concordance} If a link $L_\G$ is obtained from a link
$L$ by surgery along a connected unitrivalent graph
$\G$ with  $\ell(\G)\geq 1$ then $L_\G$ is ribbon concordant to $L$.
\end{theorem}

This result says that in the 4-dimensional setting, it is best to consider
rooted trivalent {\em trees} instead of all graphs as operators. We will thus
concentrate on trees from now on (by the STU-relation trees are sufficient
for the Vassiliev theory
$\F^v_*$ as well, see \cite{h2}). Note that grope degree (or grope class)
agrees in this case with the Vassiliev degree. In this setting the Vassiliev
invariant of degree 2, corresponding to the letter $Y$ has been  generalized
to an invariant of immersed
$2$-spheres in $4$-manifolds by Schneiderman and Teichner \cite{st}. It takes
into account the fundamental group (i.e. the edges of the letter $Y$  are
labeled by elements of the fundamental group with a holonomy relation around
each trivalent vertex) and is a second order obstruction for embedding a
single
$2$-sphere or mapping several
$2$-spheres disjointly into a $4$-manifold. It is expected that higher order
invariants of this type can be constructed for all labeled trees, modulo
antisymmetry, holonomy and IHX-relations.

Returning to knots, it turns out that a slight refinement of the theory
{\em does} allow concordance invariance. The idea is to allow surgery only
along {\em symmetric} trees, corresponding to {\em symmetric} grope
cobordism. These are related to the derived series rather than the lower
central series of the fundamental group. Symmetric trees have a new
complexity called the {\em height} $h$ (and the class $c$, defined for
any rooted tree, is given by the formula $c=2^h$).

\begin{figure}[ht]
\begin{center}
\epsfig{file=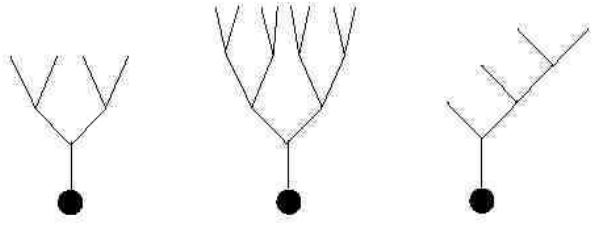,%
            height=3cm}
\end{center}
\caption{The symmetric trees of height 2 and 3, and a non-symmetric tree of
class 4}
\lbl{halffull}
\end{figure}

Recently, Cochran, Orr and Teichner defined a highly nontrivial filtration
$\F_{(h)}$ of the knot concordance group. They prove in \cite[Thm.8.11]{cot}
that two knots represent the same element in $\F_{(h)}$ if they cobound an
embedded symmetric grope of height
$\geq (h+2)$ in
$\R^3 \times [0,1]$. It is clear that a symmetric grope cobordism in 3-space
can be used to obtain such a grope cobordism in 4-space. Thus the following
consequence implies that the Casson-Gordon invariants vanish on $\F_4^{sym}$
and that the higher order von Neumann signatures of
\cite{cot} are invariants of
$\K/\F_5^{sym}$.

\begin{corollary}\lbl{4d} Define a filtration $\F_h^{sym}$ on $\K$ by
allowing symmetric trees of height
$\geq h$ as operators. Then the natural map from $\K$ to the knot concordance
group maps $\F_{h+2}^{sym}$ to the term $\F_{(h)}$ in the
Cochran-Orr-Teichner filtration.
\end{corollary}

\subsection{Open problems}
Instead of studying symmetric gropes, one can also restrict attention to any
particular grope type, parameterized by the underlying rooted tree type. The
precise definition can be found in Section~\ref{sec:equivalence}. What
follows is a summary of our low degree calculations whose proofs will be
found in \cite{ct}.

\vspace{2mm}
\noindent\begin{tabular}{|r||c|c|c|}
\hline Tree Type $\T$ &
$\K/c\T$ &
$\K/\T$ &
$\K/\T^{4}$\\
\hline\hline
\epsfig{file=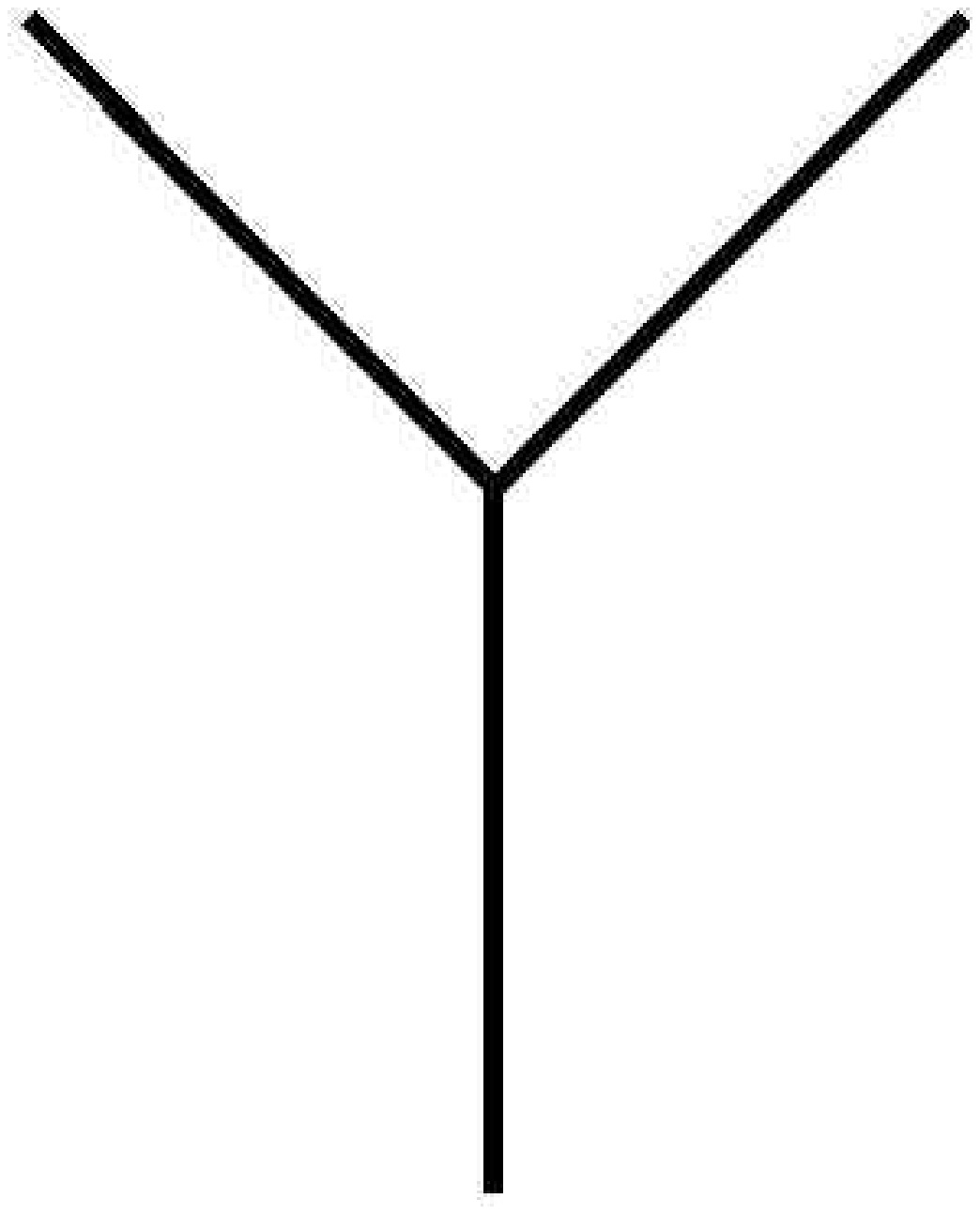,%
            height=.3cm} &
$\{0\}$&$\{0\}$&$\{0\}$\\
\hline
\epsfig{file=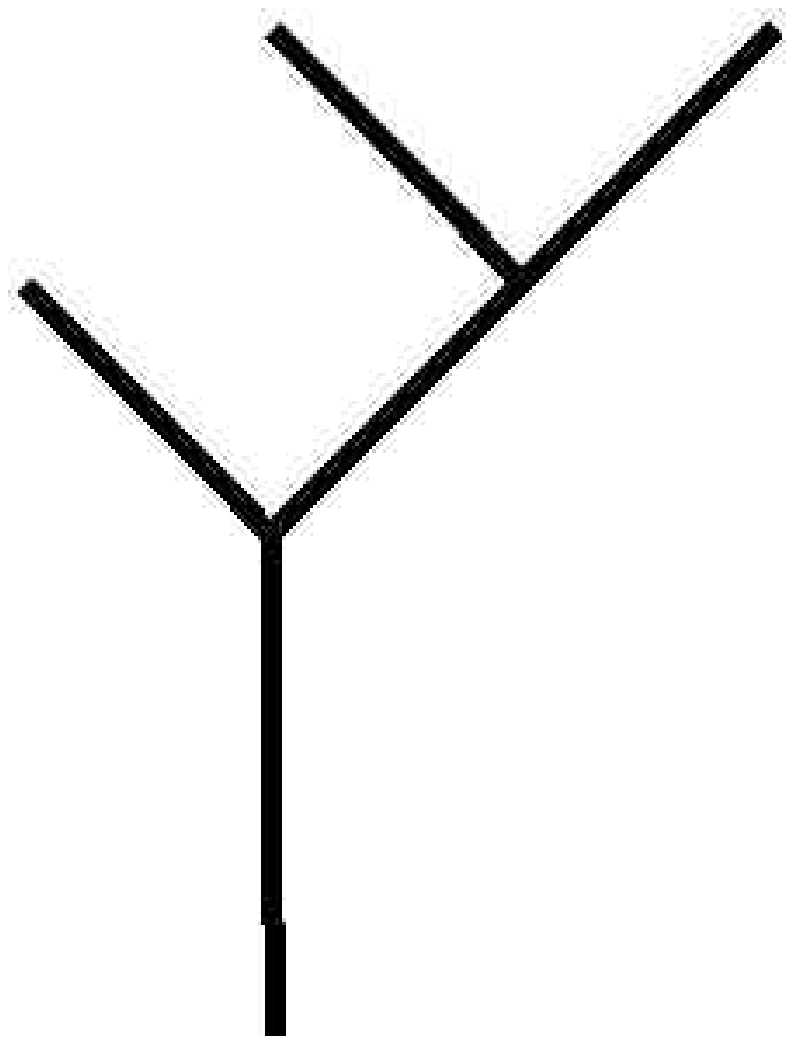,%
            height=.5cm} &
$\mathbb Z (c_2)$ &
$\mathbb Z_2 (\text{arf})$ &
$\mathbb Z_2 (\text{arf})$\\
\hline
\epsfig{file=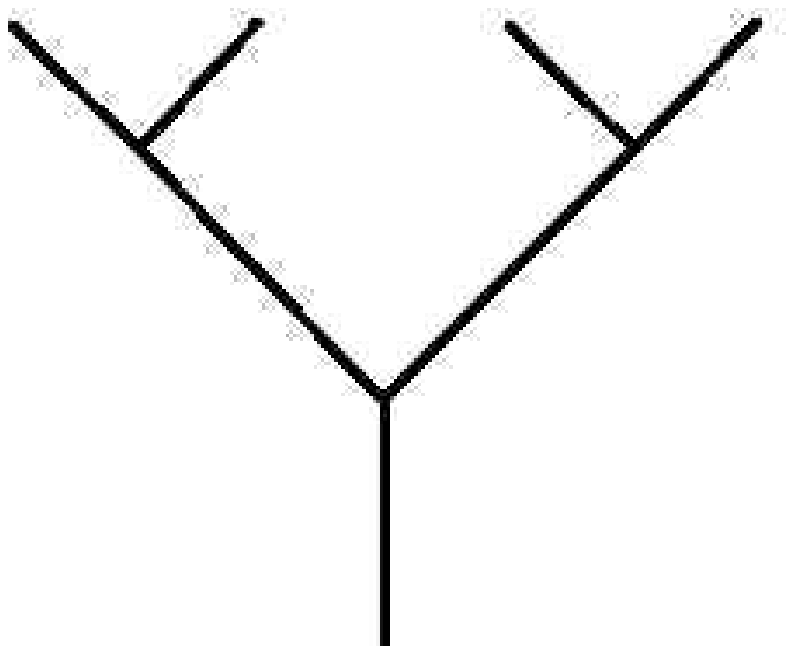,%
            height = .5cm} &
$\mathbb Z (c_3)\oplus \mathbb Z (c_2)$ &
$\mathbb Z (c_2)$ &
$\mathbb Z_2 (\text{arf})$\\
\hline
\epsfig{file=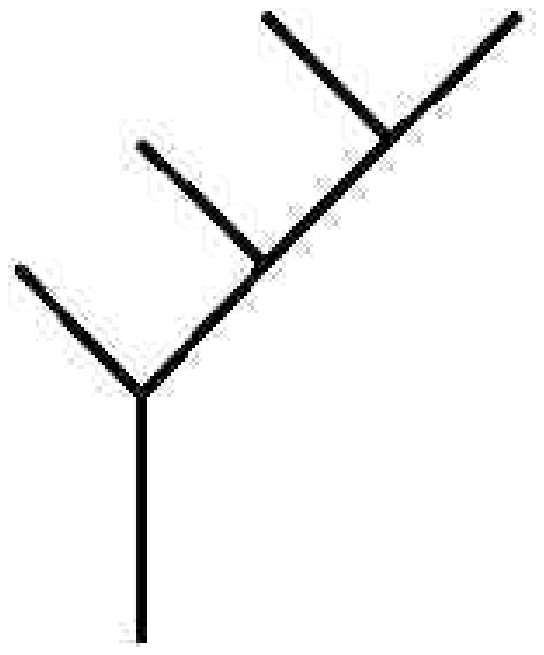,%
            height = .5cm} &
$\mathbb Z (c_3)\oplus \mathbb Z (c_2)$ &
$\mathbb Z (c_2)$ &
$\mathbb Z_2 (\text{arf})$\\
\hline
\epsfig{file=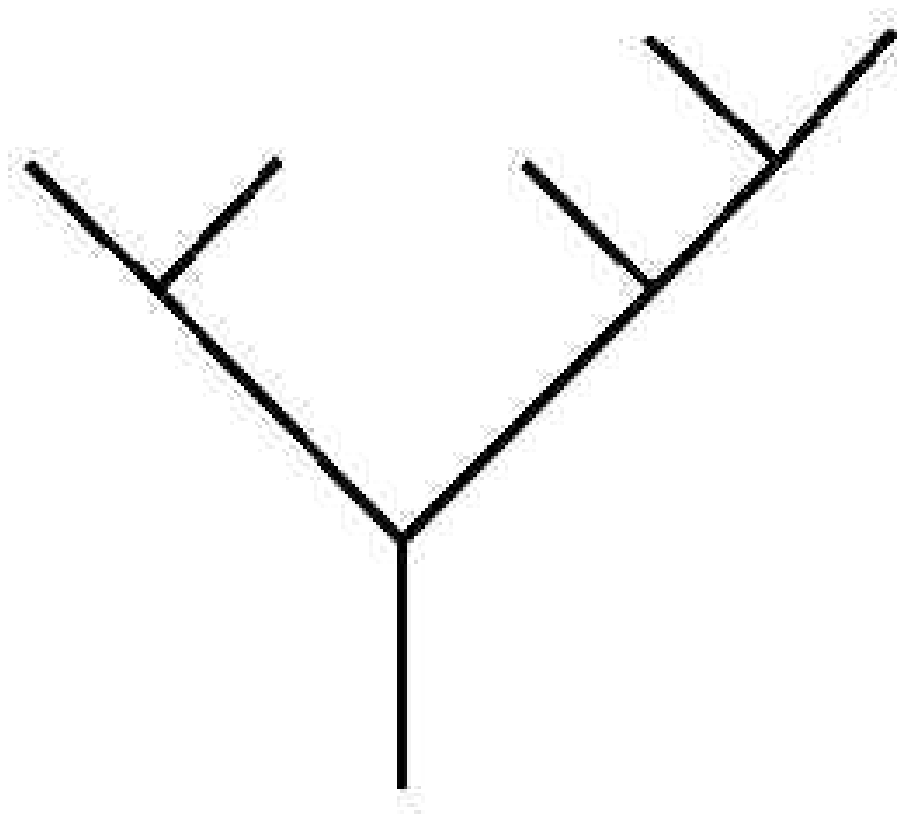,%
            height= .8cm} &
$\begin{array}{c}
\mathbb Z (c_4)\oplus \mathbb Z (c_4^\prime) \oplus\\
\mathbb Z (c_3)
\oplus \mathbb Z (c_2)
\end{array}$ &
$\mathbb Z_2 (c_3) \oplus \mathbb Z (c_2)$&
$\mathbb Z_2 (\text{arf})$\\
\hline
\epsfig{file=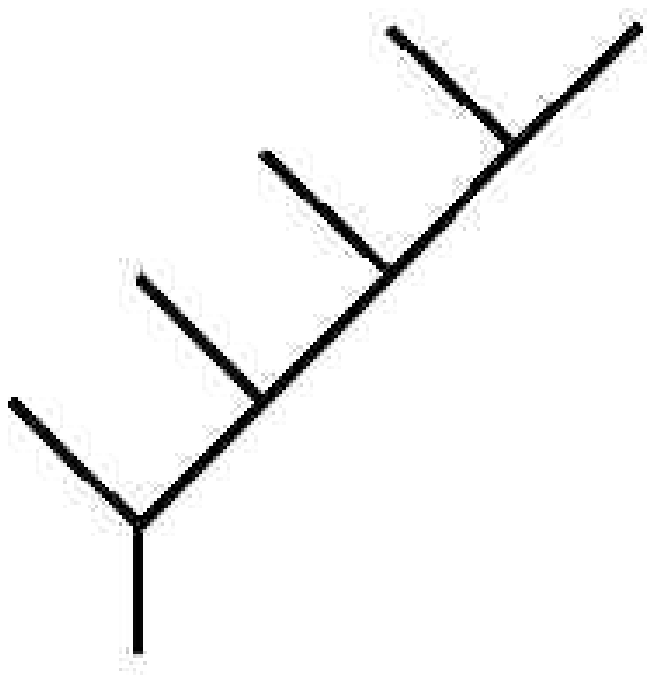,%
            height= .8cm} &
$\begin{array}{c}
\mathbb Z (c_4)\oplus \mathbb Z (c_4^\prime) \oplus\\
     \mathbb Z (c_3)
\oplus \mathbb Z (c_2)
\end{array}$ &
$\mathbb Z_2 (c_3) \oplus \mathbb Z (c_2)$&
$\mathbb Z_2 (\text{arf})$\\
\hline
\epsfig{file=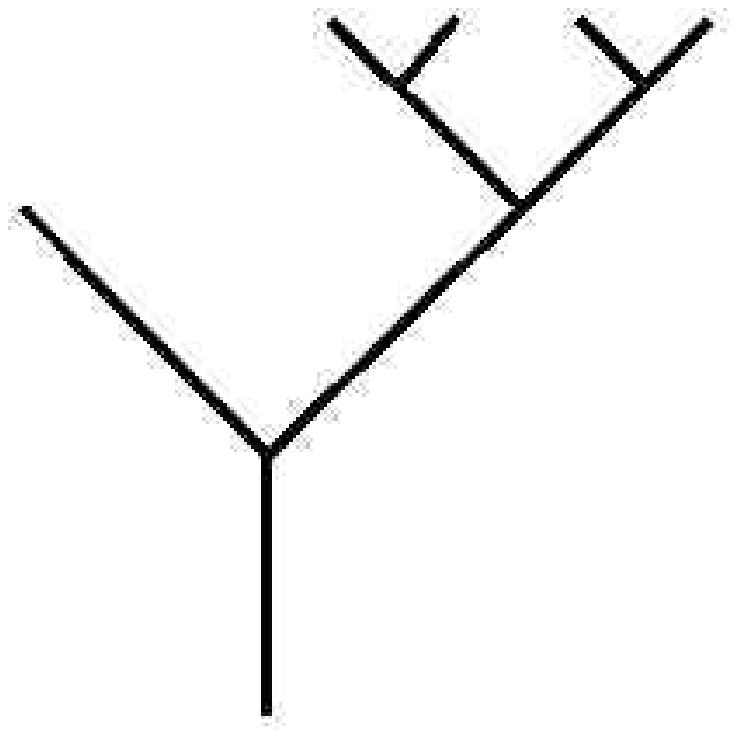,%
            height =.8cm} & {\bf ?}&
$\begin{array}{c}
\text{$S$-equivalence}\\
\text{or Bl-forms}
\end{array}
$ &
$\begin{array}{c}
\text{cobordism of}\\
\text{Bl-forms}
\end{array}
$\\
\hline
\end{tabular}
\vspace{2mm}

Here Bl-form stands for the Blanchfield form which is the equivariant linking
form on the infinite cyclic cover of the knot complement. The notation
$(c)\T$ refers to the equivalence relation given by (capped) grope cobordisms
in $3$-space using gropes of tree-type
$\T$, as explained in Section~\ref{sec:gropes}. One can also study grope
cobordism in $\R^3 \times [0,1]$ which is denoted by $\T^4$ above.

Observe that all sets in the above table are actually abelian groups (under
$\#$), except for the last row. In this case, $\K/\T^4$ is the
``groupification'' of $\K/\T$ in the sense that only the relations
$$ K + (\text{ reversed mirror image of } K)=0
$$ are added. Note that in general, this can only be true rationally because
of the occurrence of $c_3\mod 2$ in the above table. These calculations also
imply that $c_3\mod 2$ is an invariant of S-equivalence, a fact which cannot
be true rationally by \cite{mo}.

We would like to finish with the following questions and challenges for the
reader.
\begin{enumerate}
\item Find invariants of $\K/\F_h^{sym}$ for $h\geq 4$.
\item Find a good notion of grope cobordism allowing non-orientable surface
stages.
\item Can one express $4$-dimensional grope cobordism $\K/\T^4$ in terms of
algebraic operations, like the above relations, on the $3$-dimensional sets
$\K/\T$?
\item A central tool in our work is the algorithm in Theorem~\ref{refine}, which reduces
every clasper surgery to a sequence of simple clasper surgeries. It would be very useful to
implement this algorithm on a computer.
\end{enumerate}

\noindent {\em Acknowledgments:} It is a pleasure to thank Stavros
Garoufalidis and Kazuo Habiro for helpful discussions.

\tableofcontents

\section{Gropes}\lbl{sec:gropes}

\subsection{Basic definitions}\lbl{sec:gropesdefs} 
Gropes are certain
$2$-complexes formed by gluing layers of punctured surfaces together. In
our  context, a punctured surface is defined to be a closed oriented  surface
with an open disk deleted. Gropes are defined recursively using a quantity
called \emph{depth}. This differs from  the definitions in \cite{ft2} only in that
it is formally correct.

A  {\em grope} is a special pair (2-complex,circle), where the circle  is
referred to as the {\em boundary} of the grope. There is an  anomalous case
when the depth is
$1$: the unique grope of depth $1$  is the pair (circle,circle). A grope of
depth $2$ is a punctured  surface with the boundary circle specified. To form
a grope $G$ of  depth $n$, take a punctured surface, $F$, and prescribe  a
\emph{symplectic basis} $\{\alpha_i ,
\beta_j\}$. That is, 
$\alpha_i$ and $\beta_j$ are embedded curves in $F$ which represent a  basis
of $H_1 (F)$ such that the only intersections among the 
$\alpha_i$ and $\beta_j$ occur when $\alpha_i$ and
$\beta_i$ meet in  a single point. Now glue gropes of depth $<n$ along their
boundary  circles to each
$\alpha_i$ and $\beta_j$ with at least one such added  grope being of depth $n-1$. (Note
that we are allowing any added  grope to be of depth $1$, in which case we
are not really adding a  grope.)

\begin{definition} The surface $F\subset G$ is called the 
\emph{bottom stage} of the grope and its boundary is the boundary of  the
grope.
\end{definition}

\begin{definition} The \emph{tips} of  the grope are those symplectic basis
elements of the various  punctured surfaces of the grope which do not have
    gropes of depth 
$>1$ attached to them.
\end{definition}

For instance in Figure 
\ref{gropeman} there are $9$ tips. Depth was just a tool in defining  gropes.
More important is the
\emph{class} of the grope, defined  recursively as follows.

\begin{definition} The class of a depth $1$  grope is $1$. Suppose a grope
$G$ is formed by attaching the gropes  of lower depth
$\{ A_i, B_j\}$ to a symplectic basis $\{\alpha_i , 
\beta_j\}$ of the bottom stage $F$, such that $\bo A_i = 
\alpha_i$ , $\bo B_j =
\beta_j$. Then
$$\text{class}(G) := 
\underset{i}{\min} \{\text{class}(A_i) + 
\text{class}(B_i)\}.
$$
\end{definition}

\begin{figure}[ht]
\begin{center}
\epsfig{file=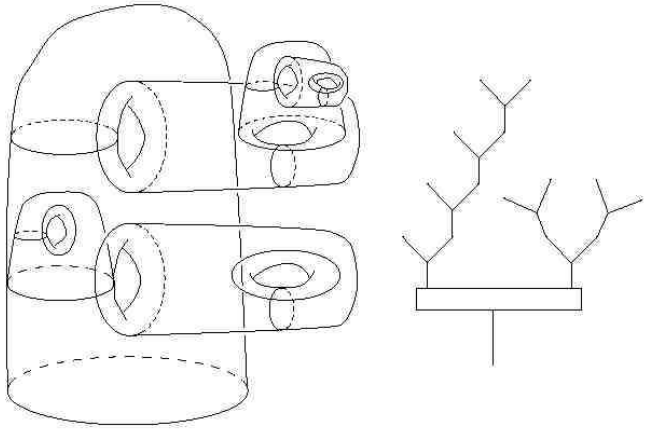,%
        height = 6cm}
\end{center}
\caption{A grope of class 4 and depth 5,  and its associated  rooted
tree-with-boxes.}
\lbl{gropeman}
\end{figure}

Associated to  every grope is a rooted tree-with-boxes. This tree is
constructed by  representing a punctured surface of genus $g$ by the
following  figure:
\begin{center}
\epsfig{file=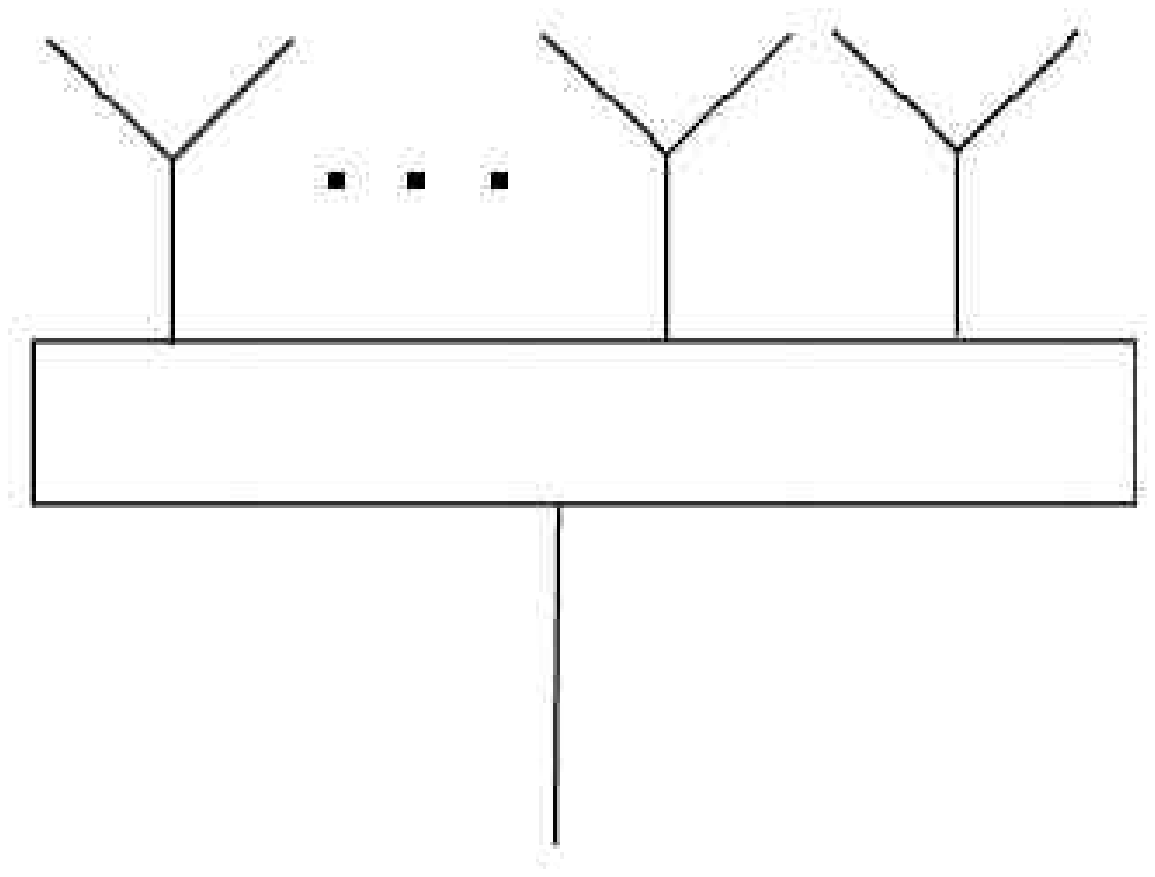,%
           height =  2cm}
\end{center} The bottom vertex is the root and it represents the  boundary of
the surface. There are
$g$ of the
\epsfig{file=wye,%
        height =.3cm} trees and the $2g$ tips of the
\epsfig{file=wye,%
        height =.3cm} 
trees represent the symplectic basis of the stage,  with dual
basis elements paired according to the \epsfig{file=wye,%
                                                height =.3cm} 
structure. Then
we glue all these trees together as  follows. If a stage $S$ is glued to a
symplectic basis element of  another stage, then identify the root vertex of the
$S$ tree, with  the tip of the other tree representing that symplectic basis
element.  Also, by convention, if a stage is genus $1$, we drop the box  and
represent that stage by a
\epsfig{file=wye,%
           height  =.3cm}.

For instance the rooted tree-with-boxes associated to the  grope in Figure
\ref{gropeman} is given on the right in that figure.  Note that depth of a
tip is the distance to the root. We will show in  Section~\ref{sec:grope
refinement} that for our purposes it is enough  to understand gropes of genus
one, i.e. gropes such that all surface  stages have genus one. These can be
represented by rooted trees  (without boxes) on which we concentrate from now
on.

A very special  class~$k$ grope is the class~$k$
\emph{half-grope} (of genus one). It  corresponds to a right-normed
commutator of length $k$ at the lower  central series level. The class $2$
half-grope tree is just a 
\epsfig{file=wye,%
           height =.3cm}. The class~$k$ half-grope  tree type is defined
recursively by adding a class
$k-1$ half grope  tree to one of the two tips of a \epsfig{file=wye,%
           height  =.3cm},  see Figure~\ref{fig:half}.

\begin{figure}[ht]
\begin{center}
\epsfig{file=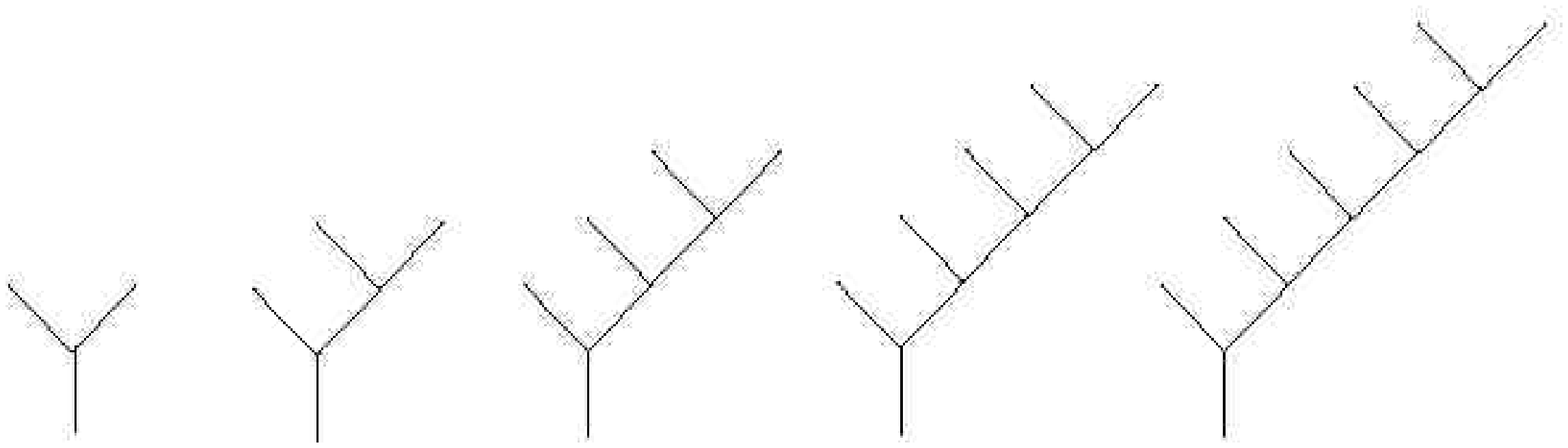,%
        height=3cm}
\end{center}
\caption{The half-gropes of class 2 to 6 and  genus 1.}
\lbl{fig:half}
\end{figure}

From these definitions, the reader should now be able to prove the following result, see
also \cite{ft2}.
\begin{proposition}
Given a continuous map $\phi:S^1\to X$, the following statement are equivalent for each
integer $k\geq 2$:
\begin{enumerate}
\item $\phi$ represents an element in $\pi_1X_k$, the $k$-th term of the lower central
series of $\pi_1X$.
\item $\phi$ extends to a continuous map of a half-grope of class $k$ into $X$.
\item $\phi$ extends to a continuous map of a grope of class $k$ into $X$.
\end{enumerate}
\end{proposition}

There are also \emph{symmetric gropes}, corresponding by a theorem just like above to the
derived series of a group, as opposed to the lower central series. A 
\epsfig{file=wye,%
           height =.3cm} represents a symmetric  grope of class
$2$. Inductively, a symmetric grope tree of  class
$2^n$ is formed by gluing symmetric gropes of class $2^{n-1}$  to the two
tips of a \epsfig{file=wye,%
           height =.3cm} as  in Figure \ref{fig:full}. A symmetric grope of
class
$2^h$ is said to  be of \emph{height} 
$h$.

\begin{figure}[ht]
\begin{center}
\epsfig{file=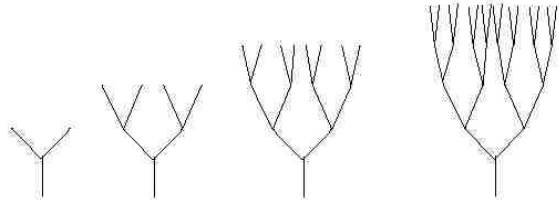,%
        height=3cm}
\end{center}
\caption{The symmetric gropes of height 1 to  4 and genus 1.}
\lbl{fig:full}
\end{figure} Sometimes we consider a  grope to be augmented with
\emph{pushing annuli}. A pushing annulus  is an annulus attached along one
boundary component to a tip of the  grope as in Figure~\ref{fig:gropecob}. It
is clear that every  embedding of a grope into 3-space can be extended to an
embedding of  the augmented  grope.
\begin{figure}[ht]
\begin{center}
\epsfig{file=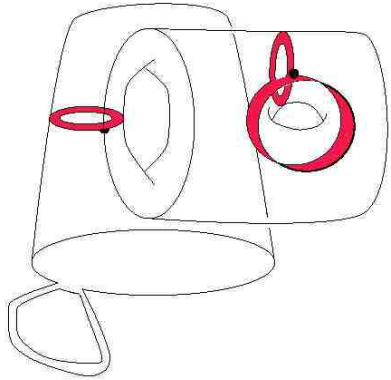,%
         height=5cm}
\end{center}
\caption{A grope cobordism with pushing  annuli.}
\lbl{fig:gropecob}
\end{figure}

\begin{definition} A {\em  capped} grope is a grope with disks (the {\em caps})
attached to all its  tips. The grope without the caps is sometimes called the
{\em body}  of the capped grope and the rooted tree type is unchanged by 
attaching caps.
\end{definition}

The (capped) gropes we have just  described have a single boundary circle, a
fact that was convenient  in the inductive definitions. But in general we
allow (capped) gropes  with an arbitrary closed 1-manifold as boundary. Such
gropes  are obtained from a grope as above by deleting open disks from the 
bottom surface stage. In particular, the relevant gropes for a {\em  grope
cobordism} between two knots will have two boundary components  as in
Figure~\ref{fig:grope}. They can also be viewed as gropes  with a single
boundary circle with an annulus attached as in  Figure~\ref{fig:gropecob}. By
definition, removing disks from the  bottom stage does not alter the
corresponding rooted tree, and adding  caps does not change the boundary of 
the grope.

\subsection{Grope cobordism of knots in $3$-manifolds}\lbl{sec:grope cobordism} 
Fix an oriented $3$-manifold $M$ and recall the basic Definition~\ref{def:grope cobordism}
from the introduction. Let $\K_M$ be the set of oriented knot types, i.e. isotopy classes of
oriented knots, in $M$.

\begin{definition} \lbl{def:M grope cobordism}
Two knot types $K_1$ and $K_2$ in $\K_M$ are {\em grope cobordant}
if there is an embedding of a grope $G$ with two boundary
components into $M$ so that the restrictions of the embedding to the two boundary
components represent the knot types $K_i$. The grope $G$ is also called a {\em  grope
cobordism} between $K_1$ and $K_2$. Only the orientation of the  bottom stage of $G$ is
relevant.
\end{definition}

It is essential to  note that the two boundary components of a grope cobordism
may link  nontrivially in $M$.
Gropes  have natural complexities associated to them. One way to make this precise is used
in our main results, Theorems~\ref{capped}  and \ref{main}:

\begin{definition} \lbl{def:class grope cobordism} 
Consider $K_1,K_2\in \K_M$ and fix an integer $c\geq 2$.
\begin{itemize}
\item[(a)] $K_1$ and $K_2$ are {\em grope cobordant of class $c$} if there
is a grope cobordism of class~$\geq  c$ between them.
\item[(b)] If there is a grope cobordism $G$ of class $c$ between $K_1$ and 
$K_2$ which extends to a map of a capped grope, such that the (interiors of
the) caps are embedded disjointly and only intersect $G$ along the bottom
stage, then $K_1$ and $K_2$ are called {\em capped} grope cobordant of class~$c$.
\end{itemize}
\end{definition}

\begin{remark}
If we don't allow the caps to intersect the grope body in  (b) then one can
do surgery on the grope (along a choice of caps) to  turn it into an annulus,
implying that $K_1$ and $K_2$ are isotopic.  Therefore, one has to somehow
weaken the notion of an embedded capped  grope. 

In dimension~$4$ one
considers {\em proper immersions} of a  capped grope \cite[\S 2.2]{fq}. This
means that the grope body is  embedded, the caps are disjoint from the body,
but the caps can  self-intersect and intersect each other. However, this
notion cannot  be useful in dimension~$3$ because of Dehn's lemma: Immersed
disks in  $3$-manifolds can usually be promoted to embedded disks, thus again  giving
an isotopy between $K_1$ and $K_2$ in our context.

This is the reason why we picked the above definition (b). Asking that the
caps only intersect the bottom stage simplifies the discussion, and is
inspired by dimension~$4$, where one can always {\em push  down}
intersections along the grope \cite[\S 2.5]{fq}. 
It turns out that in our 3-dimensional discussion the same exact statement is true: 
\end{remark}
There is another natural definition of ``capped," as suggested by the above remark,
but this turns out to be the same as the one we give: 

\begin{theorem}
Two knot types are capped grope cobordant of
class~$c$ if and only if there is a grope cobordism of class~$c$ with disjointly embedded
caps (intersecting the grope body in an arbitrary way).
\end{theorem}

 The proof of this theorem
is much more difficult than in dimension~$4$ and it requires a careful analysis of all the
steps in the proof of Theorem~\ref{capped}. We leave this proof to the interested reader.

We were so careful about {\em knot types} versus actual knots in
Definition~\ref{def:M grope cobordism} because we wanted the following Lemma to hold. Recall
that not even the relation ``two knots cobound an embedded annulus'' is an equivalence
relation on the space of knots. Therefore, one needs to work modulo isotopy all along.

\begin{lemma}\lbl{lem:equivalence relation} The relations (a) and (b) are
equivalence relations on $\K_M$.
\end{lemma}

\begin{proof}
Symmetry holds by definition. An annulus can be used to produce a grope of
arbitrary class by gluing a  trivial standard model into a puncture. Thus
annuli can be used to  demonstrate reflexivity in all cases. Transitivity
should follow from  gluing two grope cobordisms together. This can be done
ambiently in $M$ but extra care has to be taken to keep the glued grope embedded.  For
case (a) it can be seen as follows.

One proves by induction on the number of surface stages that a grope cobordism $G\subset M^3$
can be isotoped arbitrary close to a 1-dimensional complex $g\subset G$. One may assume that
this {\em spine} $g$ contains all the tips
of the grope and {\em one} boundary circle $\bo_0 G$.
Since we may use a strong deformation retraction of $G$ onto $g$, the spine $g$ (and in
particular, $\bo_0G$) is not moved during the isotopy. However, the other boundary circle
$\bo_1G$ then undergoes quite a complicated motion and ends up running parallel to all
of $g$. This means that in the following we have to be careful about introducing crossing
changes on $g$ because that might change the knot types of {\em both} boundaries of
$G$.

To prove transitivity of grope cobordism, we assume that two grope cobordisms $G,G'\subset
M$ of class $c$ are given with the knots $\bo_0G$ and $\bo_0G'$ being isotopic. After pushing
long enough towards the spines $g,g'$, we may assume that $G$ and $G'$ are disjointly
embedded. This isotopy does not change the knot types on the boundaries of the gropes, even
though it may change the 4 component link type of these boundaries (but that's irrelevant for
our purposes). We now use our assumption and start moving the knots $\bo_0G$ and $\bo_0G'$
closer to each other until they are parallel. At this point, we have to be careful not to
change the knot types of $\bo_1G$ and $\bo_1G'$, e.g. we can't just arbitrarily push $\bo_0G$
around in $M$. In fact, as pointed out above, $\bo_0G$ may {\em not} cross $g$ at all,
whereas it can cross $g'$ without changing any knot types. The same applies vice versa to
$\bo_0G'$. 

To avoid changing our knot types, we first embed an isotopy between $\bo_0G$ and $\bo_0G'$
into an ambient isotopy and run it until these knots are parallel, but
with possibly parts of
$g,g'$  still sitting in between them. Then we push $g$ across $\bo_0G'$ and $g'$ across
$\bo_0G$ until $\bo_0G$ and $\bo_0G'$ are honestly parallel in $M$. Finally, we consider tiny
thickenings of the newly positioned $g$ and $g'$ to gropes and glue them together using our
pararellism. This may require twisting the annular region around, say $\bo_0G$, so that the
gluing in fact produces an embedded grope of class~$c$ as desired. Notice that the twisting 
does not affect the isotopy class of $\bo_0G$ or $\bo_1G$.

Now we turn to case (b), i.e. transitivity of capped grope cobordism. We use the same
notation as in the previous case. In addition, we denote by $C_1,\dots,C_n$ the caps of the
grope $G$. By definition, the boundaries $\bo C_i$ are the tips of the grope 
and hence contained in the spine $g$. Hence the isotopy which pushes $G$ towards the spine
$g$ can be done relative to $\bo C_i$ and we decide to do this isotopy with {\em
all of} $C_i$ fixed. This implies that the relevant data are the disjointly embedded caps
$C_i$ (except for the usual intersection points on $\bo C_i$), together with the spine
$g$ which intersects the interiors of the caps. 

Next we implement the assumption that the
caps only intersect the {\em bottom} stage of the original grope $G$. Since $g\subset G$ this
will still be true for a tiny neighborhood of $g$ in $G$, which we now proceed to call $G$.
By general position, the intersections of this thin grope $G$ with the interiors of the caps
are thus given by short arcs which run either from $\bo_1G$ to $\bo_1G$, or from
$\bo_0G$ to $\bo_1G$. The case $\bo_0G$ to $\bo_0G$ does not occur because we chose $\bo_0G$ to be
part of the spine $g$. 

Before proceeding with the argument, we devote a paragraph to what happens if a cap were allowed to
intersect higher stages of the original grope. Then the intersections with the thin grope $G$ would
not be short arcs but rather certain unitrivalent trees which represent a normal slice through a
grope. For example, for each intersection with the second stage one would see a small H-shaped tree
in the cap, and the four univalent vertices of the H would lie on $\bo_1G$.
This can be illustrated in Figure \ref{fig:realisticgrope} (which is not capped). The second
surface stage is the big evident Seifert surface with two dual bands. An intersection of some disk
through one of these bands also picks up intersections with the grope's bottom stage, which has one
band which traces around the boundary of the second stage.
 Similarly, for
each intersection of a cap with the $r$-th stage of the grope one would see a small tree with
$(2^r-2)$ trivalent vertices and $2^r$ univalent vertices (which would lie on
$\bo_1G$). Thus the topology of these intersections distinguishes the different stages of the
grope. In the following, we shall refer to all such intersections with stages of the
grope above the bottom as ``H-shaped''. Only the bottom stage produces arcs of
intersections, and only in this case can the preferred boundary $\bo_0G$ appear in the
interior of a cap.

Now consider two capped gropes of class $c$ with caps $C_i$ respectively $C_j'$ and grope
bodies $G,G'$ which we already assume to be pushed close to the spines $g,g'$ (keeping the
caps constant).  We then do the same move as in the uncapped case, making $\bo_0G$ and
$\bo_0G'$ parallel in $M$. This can be done keeping the caps constant because $\bo_0G$
and $\bo C_i$ are disjoint parts of $g$, and similarly for $g'$. After twisting an annulus
as before, we may glue the grope bodies along the common boundary $\bo_0G$  to obtain an
embedded grope $G\cup G'$ of class~$c$. The intersection arcs of the caps with the glued up
annular regions (around $\bo_0G$ and $\bo_0G'$) now all run from $\bo_1G$ to
$\bo_1G'$, hitting the intersection $G\cap G'=\bo_0G$ once on the way. These are
intersections of the caps with the new grope's bottom stage, and hence are allowed.

We need to clean up the intersections of the caps which intersect each other and also the higher
stages of the new grope. These intersections are totally arbitrary, except the two sets
of caps are disjoint and the $C_i$ caps avoid the higher stages of $G$ and the $C'_j$ caps avoid
the higher stages of $G'$.
A consequence of the first fact is that
 there are no triple points of intersection among the caps.
After pushing little
fingers across the boundary of the caps, there are no circles of double
points, but we gain some new intersections of caps $C_i$ with a top stage of $G'$ and vice versa.
Now consider one cap
$C_i$ and recall that near its boundary  a normal slice of $G$ is H-shaped, with $\bo_1G$ on
the univalent vertices. This implies that we may push every intersection that does not contain this
knot $\bo_1G$ off $C_i$ and across
the normal slice. In particular, all intersections with
$C_j'$ can be removed this way: Every ribbon and clasp intersection can be pushed across the
boundary of
$C_i$ because only crossing changes between $\bo_1G$ and $\bo_1G'$ are introduced (and all knot
types stay the same). Doing this clean up procedure with each of the caps $C_i$, we end up
with disjointly embedded caps for $G\cup G'$, but possibly intersecting all stages of this
grope.

The next step, now that all the caps are disjoint, is to remove intersections of the caps $C_i$ with
higher stages of the grope $G'$, and vice versa. 
Suppose that a cap $C_i$ intersects higher stages of $G'_i$. It will do so along some 
unitrivalent graph, but any univalent vertices are part of
 $\bo_1G'$. Thus we may push all of these intersection out of the cap $C_i$ and across the normal
slice, introducing crossings of $\bo_1G'$ and $\bo_1G$, which do not change the isotopy class of
either. 
Similarly, higher stages of $G$ will only intersect caps $C_j'$ so that they
can be pushed off again without changing the knot types.
This leads to a capped grope cobordism of class~$c$ between $\bo_1G$ and
$\bo_1G'$ and thus transitivity is proven.
\end{proof}

\subsection{Grope refinement}\lbl{sec:grope refinement}

We will  presently refine the notion of grope cobordism by prescribing the
rooted tree type of the grope instead of just restricting its class.
However, it is technically easier to just do this for genus one gropes.
Therefore, we first discuss how to reduce to this case by presenting the 3-dimensional
version of a  technique discovered by Krushkal \cite{k} to refine gropes in 
$4$-manifolds into genus one gropes.

\begin{proposition}\lbl{reducegenus}
   Every (capped) grope  cobordism $G$ in $M$ can be realized as a sequence of
(capped) genus one grope  cobordisms $G_i$.
   Moreover, the rooted tree types of $G_i$ can be  obtained from the rooted
tree-with-boxes of $G$ by iteratively  applying the algorithm of
Figure~\ref{fig:pushgenus} to push boxes  (or genus) down to the  bottom.
\end{proposition}
\begin{proof}
\begin{figure}[ht]
\begin{center}
\epsfig{file=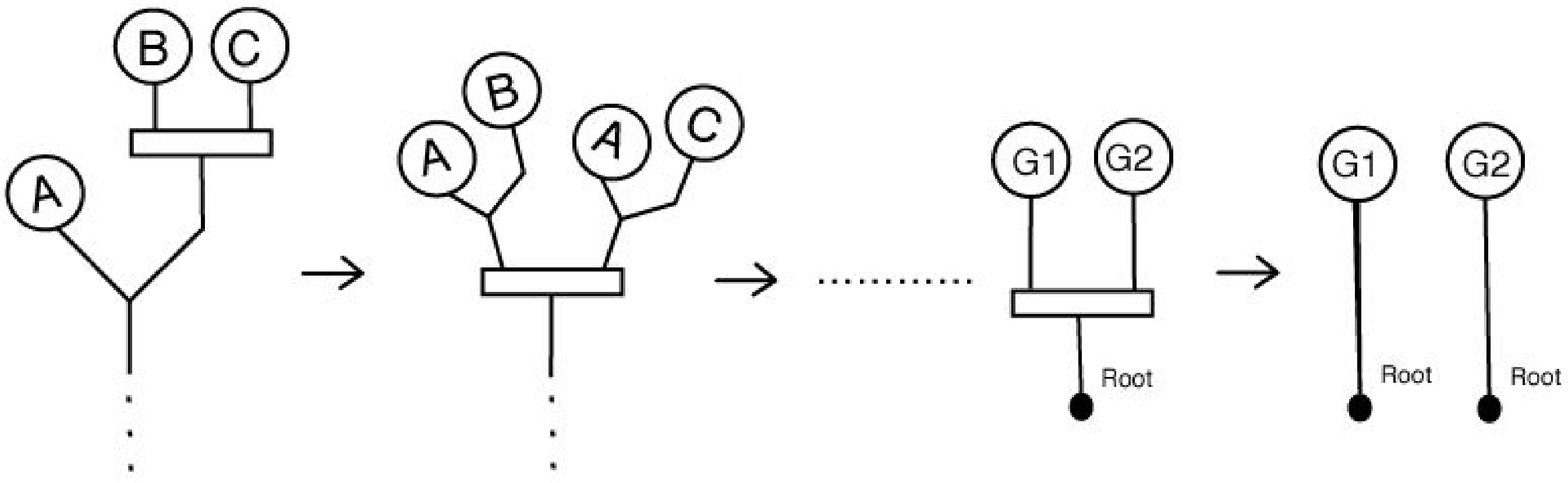,%
        height = 3cm}.
\end{center}
\caption{Pushing genus down the  grope.}
\lbl{fig:pushgenus}
\end{figure}

The way to push genus  down the grope is shown in Figure
\ref{fig:slavatrick}. It shows how  to trade genus of a stage with the
previous stage. You run an  arc from the previous stage across the current
stage in such a way as  to separate the genus. Then run a small tube along
the arc,  increasing the genus of the previous stage. The dual stage is 
depicted by $A$ in the picture. In order to make the tree type of the  grope
behave as on the left of Figure~\ref{fig:pushgenus}, we push  off a parallel
copy of $A$. (In the capped situation, $A$ will have caps, which should be included when
pushing off a parallel copy. The new caps will also only intersect the bottom
stage.)
 The parallel copies  of $A$ may intersect, a fact we have
depicted in Figure~\ref{fig:slavatrick}. 
 (In 4 dimensions, however, they do
not intersect if the grope is  framed, so there is no further  problem.)

\begin{figure}[ht]
\begin{center}
\epsfig{file=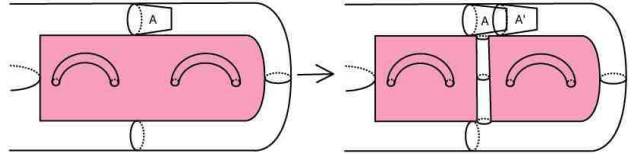,%
        height = 2.7cm}
\end{center}
\caption{Krushkal's grope  refinement.}  \lbl{fig:slavatrick}
\end{figure}

However, we can  still iteratively apply this procedure, despite the  self
intersections until all the genus is at the bottom stage. But we  can further
subdivide the resulting grope cobordism with genus $g$ at  the bottom stage
into a sequence of $g$ cobordisms with genus one at  the bottom stage, as on
the right of Figure~\ref{fig:pushgenus}.
 We  claim that each genus one grope
cobordism $G_i$ is embedded. This can  be seen schematically in Figure
\ref{fig:intpattern} which is supposed to  show that the only intersections that
arise come from parallel copies 
$A$ and
$A'$ which will eventually belong to distinct gropes $G_i$  and
$G_j$. This follows from the fact that the tree type of the  gropes only
changes as in Figure~\ref{fig:pushgenus} which implies  that at each step
parallel copies correspond to distinct branches emanating out of a box. In
the last step of the pushing down procedure, these different branches
actually become distinct gropes 
$G_i$. 

If there are caps, note that they will still only intersect the
bottom stage of the grope $G_i$ they are attached to, even though they may intersect higher
stages of $G_j,j\neq i$. 
\end{proof}

\begin{figure}[ht]
\begin{center}
\epsfig{file=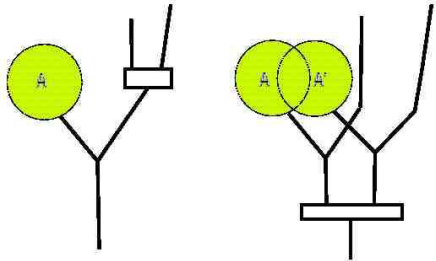,%
        height = 4cm}
\end{center}
\caption{Schematic of self-intersections  arising from grope
refinement.} \lbl{fig:intpattern}
\end{figure}

\subsection{$\T$-grope cobordism of knots in $3$-manifolds}\lbl{sec:equivalence}
We have seen in the  previous section that it is enough to consider genus one
grope  cobordisms in a $3$-manifold $M$. However, the genus of the  bottom 
surface should not be restricted to one.

\begin{definition} \lbl{type} Let $\T$ be a rooted trivalent tree. If a grope 
$G$ can be cut along the bottom surface into genus one  gropes of type $\T$
then we call $G$ a
$\T$-grope. Adding  caps to all the tips of $G$ makes it a capped
$\T$-grope.
\end{definition} 
This definition is introduced to make the following notions
of grope cobordism in $M$ into equivalence relations by composing
cobordisms. This is much more natural than taking the equivalence relation
{\em generated by} genus one grope  cobordisms of fixed tree type.
Transitivity is potentially useful for  applying 3-manifold techniques to the
study of Vassiliev  invariants.

\begin{definition}\lbl{grequivdef} Let $K_1,K_2\in\K_M$ be  oriented knot types and
$\T$ be a rooted  trivalent tree.
\begin{itemize}
\item[(a)] $K_1$ and $K_2$ are {\em 
$\T$-grope cobordant} if there is an embedding of a 
$\T$-grope into $M$ whose two boundary components represent $K_1$ and 
$K_2$.
\item[(b)] $K_1$ and $K_2$ are {\em capped
$\T$-grope  cobordant} if there is a mapping of a capped
$\T$-grope into $M$  whose boundary components are $K_1$ and $K_2$.
This  mapping is required to be an embedding except that the (disjointly embedded)  caps
are allowed to intersect the bottom stage surface of the grope.
\end{itemize}
\end{definition} 
The following result was implicitely proven in Lemma~\ref{lem:equivalence relation}:
\begin{lemma} The relations (a) and (b) are equivalence relations.
\end{lemma}

\begin{corollary} The equivalence relation  generated by {\em
genus one} (capped) $\T$-grope cobordism is exactly the same
as (capped) $\T$-grope cobordism (where the bottom stage has
arbitrary  genus).
\end{corollary}

\section{Claspers} \lbl{sec:claspers}
\subsection{Basic definitions}\lbl{sec:claspersdefs}

We recall the main notions from  Habiro's paper
\cite{h2}, making an attempt to only introduce the  notions relevant to grope
cobordism and the relation to finite type  invariants. In particular, we
completely avoid all the boxes in  claspers since we can always reduce to
this case.

A  clasper is a  compact connected surface made out of the following 
constituents:
\begin{itemize}
\item {\em edges} are bands that  connect the other two constituents,
\item {\em nodes} are disks with  three incident edges,
\item {\em leaves} are annuli with one incident  edge.
\end{itemize}
\begin{figure}[ht]
\begin{center}
\epsfig{file=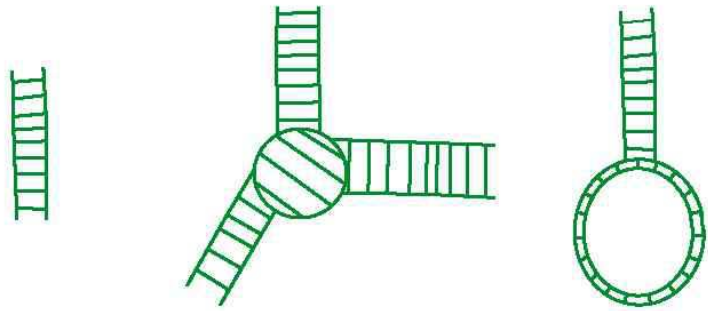,%
        height = 5cm}
\end{center}
\caption{An edge, a node and a  leaf.}
\end{figure}

Thus a clasper collapses to a unitrivalent graph  such that the nodes become
one type of trivalent vertex and each leaf  has exactly one trivalent vertex
of a second type. However, it is  common to think of this second type as a
univalent vertex (ignoring  the leaves momentarily) and to only consider
those vertices as  trivalent that come from nodes. If $\G$ is the underlying 
unitrivalent graph of a clasper (again ignoring the leaves), then we  call it
a $\G$-clasper, and we call $\G$ the {\em type} of the  clasper. A {\em tree
clasper} is a clasper whose type is a  tree.

\begin{figure}[ht]
\begin{center}
\epsfig{file=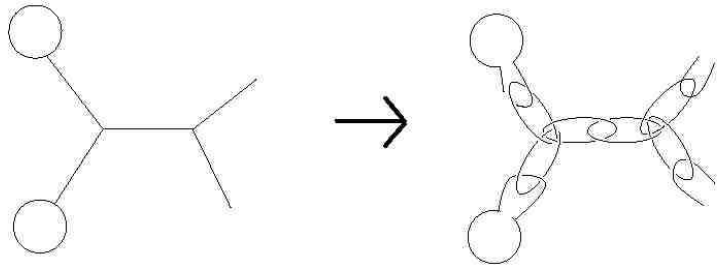,%
        height=4cm}
\end{center}
\caption{Associating a link to a  clasper.}
\lbl{claspertolink}
\end{figure}

Assume a clasper $C$ is  embedded in a $3$-manifold $M$. Then one can
associate to it a framed  link $L_C$ in $M$ by replacing  each edge by the
(positive) Hopf-link  and each node by a $0$-framed (positive) Borromean
rings, see  Figure~\ref{claspertolink}. The framing (slope along which to
attach  a $2$ handle) of each link component associated to a leaf is 
determined in the obvious way by the framing of the leaf.
 There, and  in most figures to follow, only the spine of the clasper is drawn
and  the blackboard framing is used to thicken it to a surface. Two 
thickenings differ by twistings of the bands and annuli, and also by 
reordering the three edges incident to a node. Note that a 
$0$-framing is well defined for components that lie in small balls,  usually
the  neighborhoods of a trivalent vertex or edge.

If one of the leaves of a clasper $C$ bounds a disk into $M\smallsetminus 
C$, we call it a {\em cap} because of the relation with gropes explained below. In the
presence of a cap, surgery on the framed link
$L_C$ does not change the ambient  3-manifold $M$. This implies that if $C$ lies in the
complement of a  knot $K$, then surgery on $L_C$ gives a new knot $K_C$ in the same
manifold $M$, the  {\em surgery of $K$ along $C$}. Figure~\ref{fig:intro8} shows how one  can
obtain a Figure~8 knot as surgery on the unknot along a  Y-clasper.

\begin{definition}\lbl{def:caps for claspers} 
A clasper $C$ is called {\em capped} if the leaves bound disjoint disks (the
caps) into $M \smallsetminus C$. If it happes that only some of the leaves of $C$ bound disks
into $M \smallsetminus C$ then we only call those disks {\em caps} if they are embedded
disjointly.
\end{definition}
The following notions for claspers all depend not only on the position in $M$ but also on
the relative position with respect to a knot $K$.
\begin{definition}\lbl{def:claspers} Let $C$ be a  clasper in the
complement of a knot $K\subset  M^3$.
\begin{itemize}
\item $C$ is a {\em rooted} clasper if one leaf has a cap which intersects $K$  transversely
in a single point. In particular, the surgery $K_C$ is  defined as a knot in $M$. The
particular leaf becomes also the root  of the underlying type of the clasper.
\item Conversely, if one has  given a {\em rooted unitrivalent} graph $\G$,
then a $\G$-clasper is  a rooted clasper of type $\G$.
\item If $\G$ is a  rooted unitrivalent graph then a capped 
$\G$-clasper is a capped clasper of type $\G$ such that the cap
corresponding to the root intersects the knot
$K$ transversaly  in a single point.
\item $C$ is a {\em simple} clasper if it is  capped such that each cap
intersects the knot transversaly in a  single point.
\item There are several degrees associated to claspers.  By definition, these
are the degrees of the underlying type (which  replaces the leaves by
univalent vertices). We have mentioned three  different possibilities in the
introduction, the Vassiliev, loop and  grope degrees.
\item For any such degree $\deg$, the equivalence  relation on $\K_M$ defined
by
$\F^{\deg}_k$ in the introduction is  generated by simple clasper surgeries
of degree $\deg\geq  k$.
\end{itemize}

\end{definition}
\begin{remark} The notions of  rooted and capped claspers are new and replace
notions  like admissible, strict and special in \cite{h2}. We feel that 
descriptive names are very  important.
\end{remark}

The surgery on unitrivalent graphs described in the introduction is by definition given by
clasper surgery on {\em the simple} clasper defined by the graph. Thus simple clasper
surgeries define the relevant quotients of $\K$ defined in the introduction and used in our
main Theorems~\ref{capped} and \ref{main}.

\begin{figure}[ht]
\begin{center}
\epsfig{file=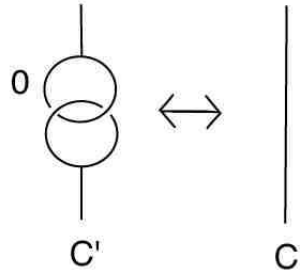,%
        height=4cm}
\end{center}
\caption{First Morse  cancellation.}
\lbl{fig:morse1}
\end{figure}

There are many  identities among claspers, perhaps the most basic of which is 
as follows. Let the clasper $C^\prime$ be obtained from $C$ by  cutting an
edge and inserting a Hopf-linked pair of tips as in  Figure~\ref{fig:morse1}.
Then surgery on
$C$ is equivalent to surgery  on $C^\prime$.  This follows from standard Kirby
calculus, or more  precisely from Morse canceling the Hopf-pair viewed as a
1-handle and  a 2-handle in the
$4$-dimensional world. 

A second often used Morse  cancellation occurs if one thinks of one of the
three Borromean rings  as a 1-handle and cancels it with a 2-handle coming
from an adjacent  leaf as in  Figure~\ref{fig:morse2}.
\begin{figure}[ht]
\begin{center}
\epsfig{file=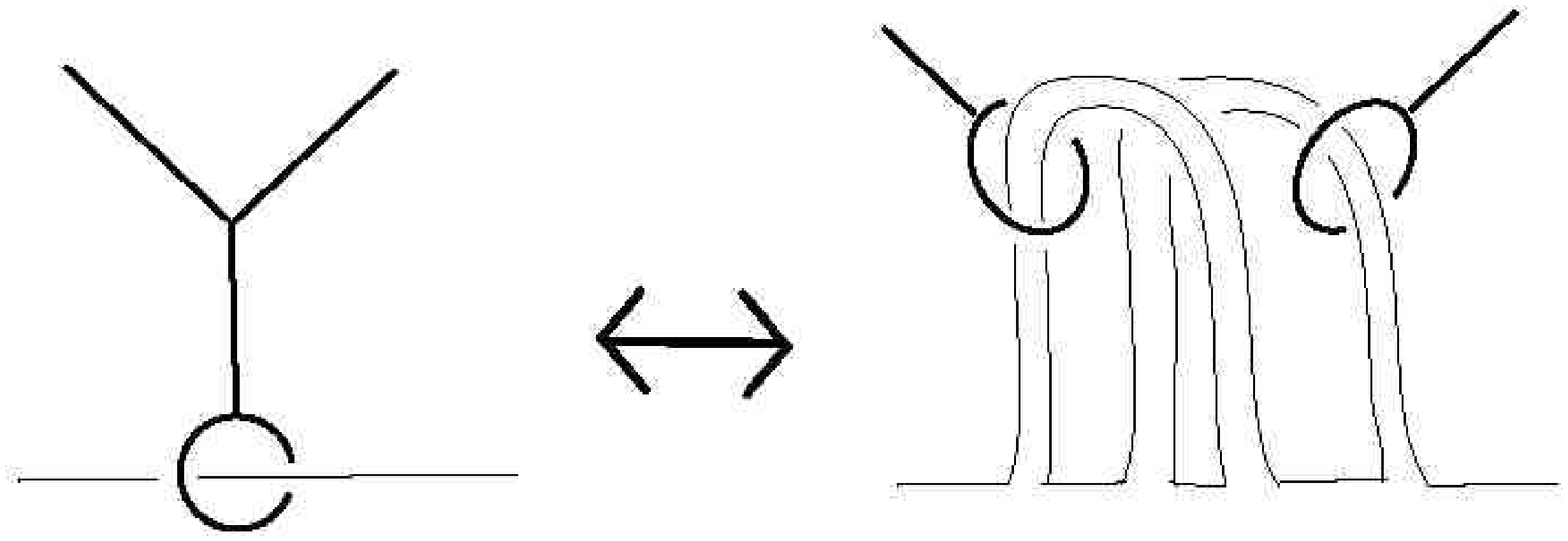,%
        height=4cm}
\end{center}
\caption{Second Morse  cancellation.}
\lbl{fig:morse2}
\end{figure}

\subsection{Claspers  and gropes}\lbl{sec:gropes and claspers}

In this section we show  that a $3$-dimensional grope cobordism of genus one
is the same as a  rooted tree clasper surgery. The rooted tree type of the
clasper is  the same as the rooted tree type of the grope. We first outline 
the construction of a clasper, given a grope cobordism, and  subsequently give
the reverse  construction.

\begin{theorem}\lbl{prevthm} Let $\T$ be a rooted  trivalent tree. Then a
$\T$-grope cobordism of genus one can be  realized by a
$\T$-clasper surgery, supported in a regular  neighborhood of the  grope.
\end{theorem}

\begin{remarks}\hspace{1em}

\begin{itemize}
\item  The clasper we obtain from the grope is not unique.  This
indeterminacy leads to a set of identities on claspers.
\item  This theorem could be strengthened to give a correspondence  between
gropes with genus and claspers with boxes, but for clarity we  do not consider
this greater  generality.
\end{itemize}
\end{remarks}

\noindent Theorem 
\ref{prevthm} will follow from the following  relative  version.

\addtocounter{thm}{-1}
\begin{thm}\lbl{string} Let $H$ be  an oriented $3$ manifold with two
distinguished points $x_0$ and $x_1$  on its nonempty boundary. Let $\alpha$ and
$\tilde{\alpha}$ be two properly  embedded arcs in $H$, with disjoint
interiors, running from $x_0$ to 
$x_1$. Suppose $\alpha
\cup \tilde{\alpha}$ bounds a $\T$-grope in 
$H$. Then there is a $\T$-clasper
$C$ embedded in $H\backslash 
\alpha$, with root a meridian to $\alpha$, such that $\alpha_C$ is isotopic to
$\tilde{\alpha}$ rel  boundary.
\end{thm}

To see that this implies Theorem~\ref{prevthm},  recall from
Figure~\ref{fig:gropecob} that a grope cobordism between  knots $K$ and
$\tilde{K}$ can be thought of as a grope $G'$ with one  boundary component,
band summed with an annulus with core (say) 
$\tilde{K}$. Consider the handlebody $H$ which is a regular  neighborhood of
$G'$. Then $K$ intersects
$H$ in an arc $\alpha$ and  the boundary $\bo H$ hits the  cobordism
along an arc 
$\tilde{\alpha}$. Together $\tilde{\alpha}\cup \alpha$ bound the  grope
$G'$ and hence there is a $\T$-clasper $C$ in $H$ which takes 
$\alpha$ to
$\tilde{\alpha}$ rel boundary. In a regular neighborhood  of the original
cobordism, $C$ therefore takes $K$ to a parallel copy  of $\tilde{K}$:
\begin{center}
\epsfig{file=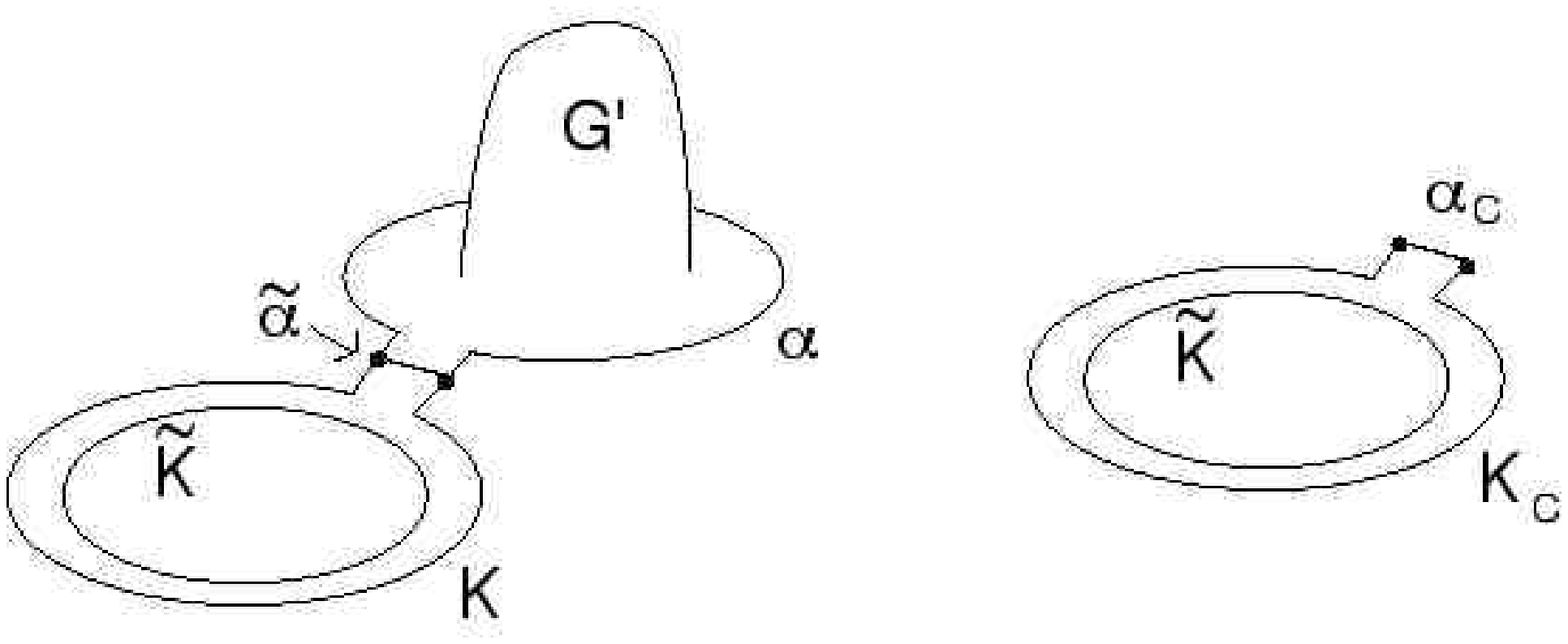,%
        height=3cm}
\end{center}

\begin{proof}[Proof of Theorem 
\ref{string}]

\noindent{\bf Construction of the unframed  clasper}

Assume the grope is augmented with pushing annuli. Then  each surface stage of
the grope has two surfaces which attach to it,  and these are either pushing
annuli or higher surface stages of the  grope. In order to simplify
terminology, refer to both these types of  surface as \emph{higher  surfaces}.

\begin{figure}[ht]
\begin{center}
\epsfig{file=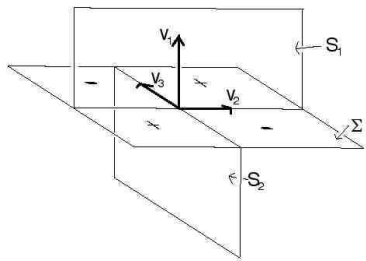,%
        height=4cm}
\end{center}
\caption{Positive  quadrants.}  \lbl{fig:posquadrant}
\end{figure}

Let $\Sigma$ be a  surface stage of the embedded grope, with higher surfaces
$S_1$ and 
$S_2$ attaching to it. Then $S_1\cap S_2$ is a point $s_0$, and in  a
neighborhood of this point $s_0$, $S_1\cup S_2$ divides $\Sigma$  into four
quadrants. We distinguish two of these as follows. Let 
$(v_1,v_2,v_3)$ be an ordered basis of the tangent space 
$T_{s_0}M$ constructed as follows. Let $v_1$ be transverse to 
$\Sigma$ and pointing into
$S_1$. Choose $v_2$ tangent to
$\Sigma\cap  S_1$. Choose $v_3$ tangent to $\Sigma\cap S_2$ in such a way 
that
$v_1\wedge v_2\wedge v_3$ is a positive orientation of $\R^3$.  The two
quadrants lying between $v_2$ and $v_3$ and between $-v_2$  and $-v_3$ are
called \emph{positive quadrants}, see  Figure~\ref{fig:posquadrant}. There
were two choices in selecting 
$v_1,v_2,v_3$, namely which surface is called $S_1$ ($v_1$ versus 
$-v_1$) and which direction of $\Sigma\cap S_1$ the vector 
$v_2$ points along ($v_2$ versus $-v_2$.) Changing $v_2$ to $-v_2$  will also
change
$v_3$ to $-v_3$ in order to preserve the orientation 
$v_1\wedge v_2\wedge v_3$. Therefore, the positive quadrants do not  change.
If one changes
$v_1$ to $-v_1$, then the role of $v_2$ and 
$v_3$ is reversed. But
$-v_1\wedge v_3\wedge v_2$ is still positive,  and hence the positive
quadrants are those between $v_3$ and $v_2$,  as before.

We are now ready to define the unframed clasper $C^u$ in 
$H\backslash
\alpha$. The leaves include those ends of the pushing  annuli which are not
attached to anything. (These are the \emph{tip  leaves}.) There is one more
leaf which is a meridian to 
$\alpha$.(This is the \emph{root leaf.}) This leaf punctures the  bottom
stage of the grope in a single point. Every surface stage  contains a node of
$C^u$ where the higher surfaces intersect. Hence  each pushing annulus has a
node on its boundary. This is connected by  an embedded arc in the annulus to
the tip leaf at the other end. Each  surface stage except the bottom stage
contains two nodes: one on  the boundary and one in the interior. Connect
these by an embedded  arc in the surface stage whose interior misses the
attaching regions  for the higher surfaces, and such that it emanates from
the interior  node in a positive quadrant. Finally connect the node on
the  bottom stage to the intersection of the root leaf with the stage by  an
embedded arc whose interior avoids the attaching regions for the  higher
surfaces, and which emanates from the node in a positive  quadrant.

Figure~\ref{fig:grope2clasper} shows the construction for  a grope of  class
$3$.
\begin{figure}[ht]\begin{center}
\epsfig{file=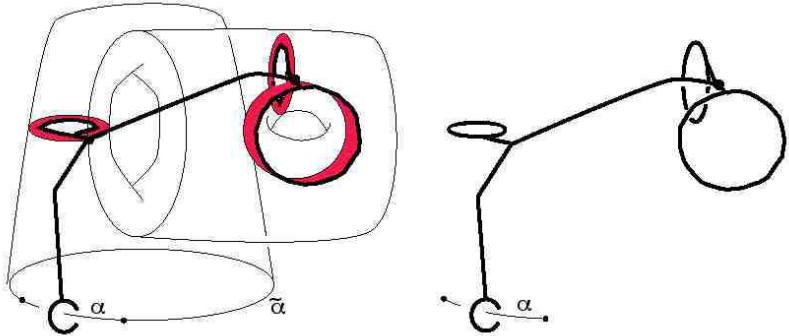,%
        height=5cm}
\end{center}
\caption{Associating an unframed clasper to  a
grope.}  \lbl{fig:grope2clasper}
\end{figure}

\noindent {\bf  Figuring out the framing}

The tip leaves of the clasper have obvious  framings along the annuli they are
contained in. Similarly each edge  has an obvious framing as a subset of  a
surface.
\begin{figure}[ht]\begin{center}
\epsfig{file=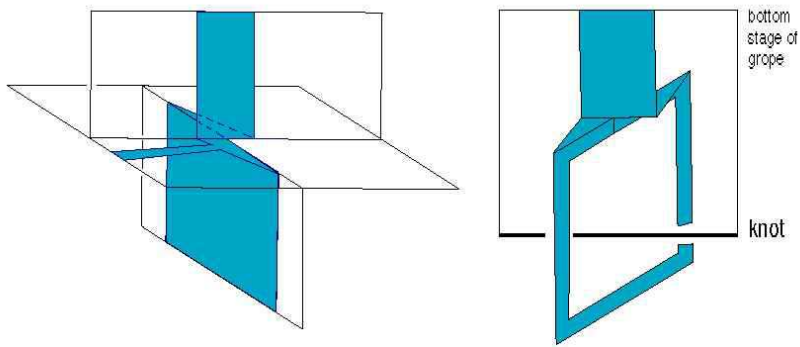,%
         height = 5cm}
\end{center}
\caption{Extending the framing to nodes  and to the root leaf.}
\lbl{fig:rootframe}
\end{figure} 

Framing a node  is depicted in Figure \ref{fig:rootframe}. Notice
that the edge on the  surface stage is approaching via a positive quadrant.
We glue  together the perpendicular framings of the two edges associated to 
the higher surfaces with two triangles inside the positive quadrants.  The
framing of the approaching edge is naturally glued to one of  these triangles.

We can frame the root leaf using the meridional  disk it bounds. This needs to
be glued to the perpendicular framing  of the incident edge.
  This is shown in  Figure~\ref{fig:rootframe},  where we again use two triangles
to glue up different parts of the  clasper. Notice that this is the only
place at which the clasper is  not a subset of the grope but the triangles
are defined as in the discussion of positive quadrants. 

\noindent {\bf Proving  that this works}

  We proceed by induction on the number of surface  stages, the base case
being a surface of genus one. Let $\Sigma$ be  the base surface, $\Sigma^a$
the augmented surface and
$C$ the  clasper we just  constructed.

\begin{figure}[ht]\begin{center}
\epsfig{file=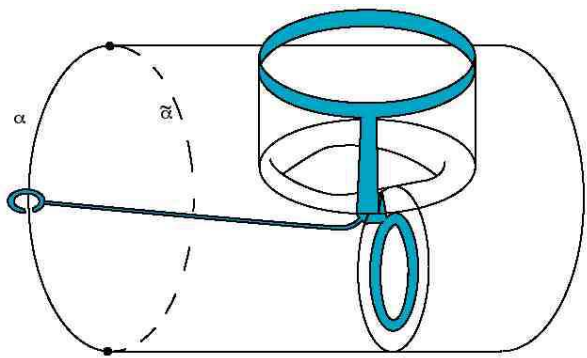,%
        height = 5cm}
\end{center}
\caption{A standard model of $(\Sigma^a,  C)$.}
\lbl{fig:stdpicture}
\end{figure}

\begin{lemma} The pair 
$(\Sigma^a , C)$ in $H$ can be realized as the restriction of  an orientation
preserving embedding into $H$ of the genus two 
handlebody, which is a regular neighborhood of the standard picture  given in 
Figure~\ref{fig:stdpicture}.
\end{lemma}
\begin{proof}
\begin{figure}[ht]\begin{center}
\epsfig{file=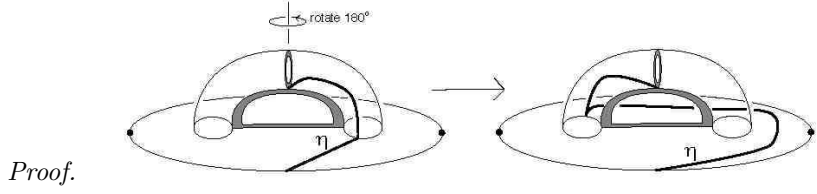,%
        height = 2.7cm}
\caption{Exchanging positive  quadrants.}
\lbl{posswitch}
\end{center}
\end{figure}

By definition $\Sigma^a$  is an embedding of the given picture, ignoring the
clasper $C$. We  precompose this embedding with a suitable orientation
preserving  automorphism of the regular neighborhood which fixes
$\bo\Sigma$ pointwise and $\Sigma^a$ setwise. Clearly the edges on  the pushing
annuli can be straightened out by twists supported in the  annuli's
interiors, and these twists extend to the regular  neighborhood. Hence it
suffices to straighten out the edge $\eta$  which runs along $\Sigma$. Let the
annuli be called $S_1$ and
$S_2$.  The interior of $\eta$ lies in the (open) annulus 
$\Sigma\backslash (\bo S_1\cup\bo S_2\cup\bo\Sigma)$.  It can
therefore be straightened via Dehn twists. It also can  approach $\bo
S_1\cup\bo S_2$ in two ways: by the two  positive quadrants. There is an
automorphism of
$\Sigma^a$ rel 
$\bo\Sigma$ taking one quadrant to the other. This is depicted  in
Figure 
\ref{posswitch}. \end{proof}

\begin{figure}[ht]\begin{center}
\epsfig{file=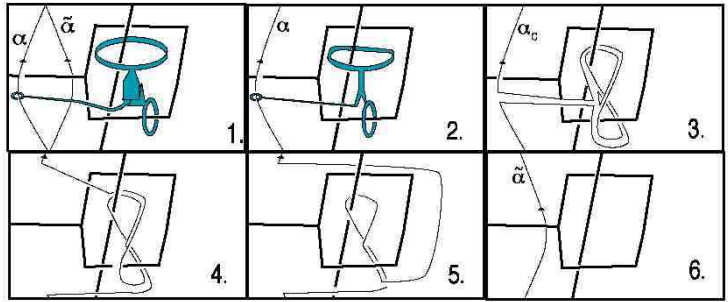,%
        height = 4.5cm}
\end{center}
\caption{The base  case.}
\lbl{fig:checkbasecase}
\end{figure}

Because of this lemma, it  suffices to check that $\alpha_C = \tilde{\alpha}$
in the standard  model of Figure~\ref{fig:stdpicture}. (We need the embedding to
preserve  orientations because an orientation is required to associate  a
well-defined link to the clasper.)

The standard model is redrawn  in Figure~\ref{fig:checkbasecase}, with heavy lines
deleted from the  ambient 3-ball to make it a regular neighborhood of
$\Sigma$. The  clasper is cleaned up a little bit in the second frame, and
then the  second Morse cancellation from Figure~\ref{fig:morse2} is used to 
produce
$\alpha_C$ in the third frame. Finally, an isotopy moves 
$\alpha_C$ to the knot
$\tilde{\alpha}$ as shown in the remaining  frames.

Now for the inductive step. This follows from  Figure~\ref{fig:inductive}.
Pictured is a top stage of the grope and part  of the clasper $C$ we 
constructed. In frame 2 we have broken the  edge of the clasper that lies on
the top surface into two  claspers
$C_T$ and $C_B$. This is the first Morse cancellation  from
Figure~\ref{fig:morse1} and gives $\alpha_C = 
\left(\alpha_{C_B}\right)_{C_T}$. By induction we know that the  clasper
surgery $C_T$ has the pictured effect on
$C_B$ since the  indicated section, $\beta$ of the leaf of $C_B$ cobounds
the surface  stage corresponding to the clasper
$C_T$ with the pictured arc 
$\tilde{\beta}$. This gives rise to a new clasper
$C^\prime  =
\left({C_B}\right)_{C_T}$ which corresponds to the grope which is  gotten by
forgetting about the indicated surface stage.
$\alpha$ and 
$\tilde{\alpha}$ still bound this new grope, and by  induction
$\tilde{\alpha} = \alpha_{C^\prime}$ which we saw is equal  to 
$\alpha_C$. \end{proof}

\begin{figure}[ht]\begin{center}
\epsfig{file=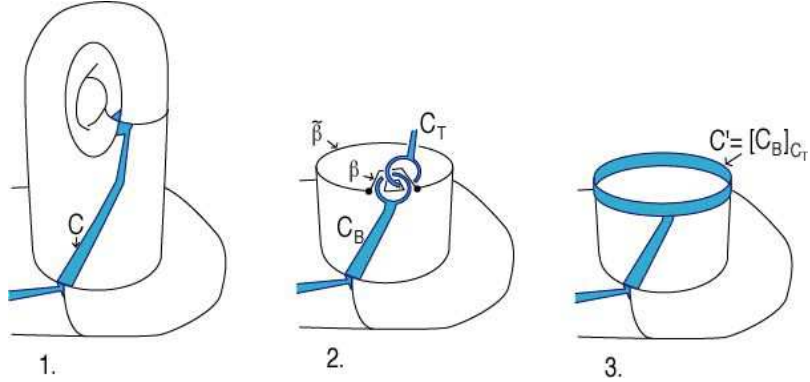,%
        height = 5cm}
\end{center}
\caption{The inductive  step.}
\lbl{fig:inductive}
\end{figure} 
We next come to the converse of  Theorem~\ref{prevthm}.
\begin{theorem}\lbl{b} Let $\T$ be a rooted  trivalent tree. Then every
$\T$-clasper surgery is realized by  a
$\T$-grope cobordism of genus one, with the grope being in a  regular
neighborhood of the clasper and knot.
\end{theorem}

As  before, it will be more convenient to prove a relative version, but 
first we introduce some notation.

\begin{definition} If $C$ is a  clasper in a 3-manifold $M$, let $M_C$ denote
the  three-manifold which is obtained by surgery on 
$C$.
\end{definition}

\addtocounter{thm}{-1}
\begin{thm}\lbl{bprime} Let 
$N$ be a regular neighborhood of a $\T$-clasper $C$. A meridian  on
$\bo N$ of the root leaf bounds a properly embedded 
$\T$-grope in $N_C$.
\end{thm}

To see that Theorem~\ref{bprime}  implies Theorem \ref{b}, suppose a
$\T$-clasper $C$ has a root leaf  on the knot $K$. Let $\tilde{K}$, be the
knot in
$M\backslash C$  where the intersection with the root leaf's disk has been
removed by  a small perturbation which pushes $K$ off that disk. Then $K$  and
$\tilde{K}$ differ by a meridian of the root leaf and hence  cobound a
$\T$-grope in $M_C$ by Theorem~\ref{bprime}. That is $K_C$  and $\tilde{K}_C$
cobound a
$\T$-grope in $M$. But $\tilde{K}_C= 
\tilde{K}=K$ in $M$, since $C$ has a disk leaf that doesn't hit 
$\tilde{K}$.

By expanding edges of claspers into Hopf-linked pairs  of leaves, Theorem
\ref{bprime} is easily seen to follow from the  following proposition.
\begin{proposition}\lbl{propproof} Let $C$  be the unique Vassiliev degree
$2$ clasper, i.e the letter Y. Let 
$N$ be a regular neighborhood of $C$. Then a meridian 
$\alpha\subset\bo N$ to any leaf bounds a properly embedded  genus one
surface in $N_C$. This surface can be augmented with two  pushing annuli
which extend to $\bo N$ as parallel copies of the  other two  leaves.
\end{proposition}

\begin{proof}

\begin{figure}[ht]
\begin{center}
\epsfig{file=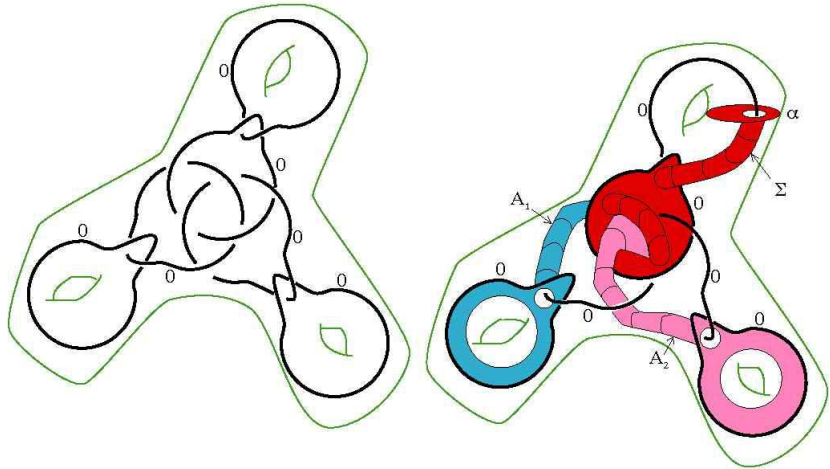,%
        height = 6.5cm}
\end{center}
\caption{The proof of proposition 
\ref{propproof}.}
\lbl{fig:clasperexpand}
\end{figure}

We have drawn 
$N$ in Figure \ref{fig:clasperexpand}, and replaced the clasper by
$0$-framed  surgery on the associated link. The curve $\alpha$ bounds the 
genus one surface
$\Sigma$. Note that part of $\Sigma$ travels over  an attached $2$ handle. Two
dual curves on $\Sigma$ each cobound an  annulus with a parallel copy of the
two lower leaves. These annuli  are denoted $A_1$ and $A_2$, and each also
runs over an attached $2$  handle. \end{proof}

\subsection{Geometric IHX and  half-gropes}\lbl{sec:IHX}

In this section we answer the question  whether grope cobordism is generated
by half-gropes, just like the  lower central series is generated by right
normed commutators. Only  for this purpose do we use concepts developed in
\cite{h2}, which  have not been covered in this paper. Denote by $\h_k$ the
rooted tree  type that corresponds to a genus one half-grope of class~$k$, as
in  Figure~\ref{fig:half}.

\begin{theorem}\lbl{half-gropes} Let $K_1$, 
$K_2$ be oriented knots in a $3$-manifold 
$M$.
\begin{itemize}
\item[(a)] $K_1$ and $K_2$ are grope cobordant  of class $k$ if and only if
there is an $\h_k$-grope cobordism  between
$K_1$ and $K_2$.
\item[(b)]
$K_1$ and $K_2$ are capped grope  cobordant of class $k$ if and only if there
is a capped $\h_k$-grope  cobordism between $K_1$ and $K_2$.
\end{itemize}
\end{theorem} The  proof of this result uses a very nice unpublished result
of Habiro,  which is a geometric realization of the IHX-relation for capped
tree  claspers.
\begin{figure}[ht]
\begin{center}
\epsfig{file=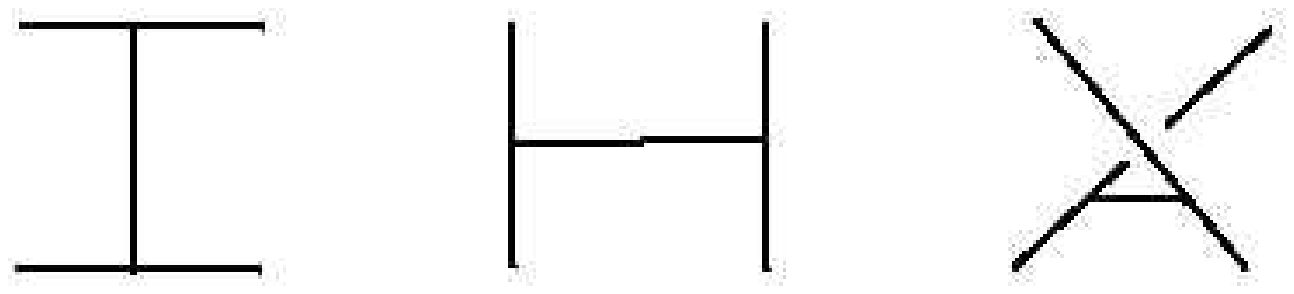,%
        height=2cm}
\end{center}
\caption{IHX.}
\lbl{fig:geomihx}
\end{figure}

\begin{theorem}  [Habiro]\lbl{Habiro} Let $I,H$ and $X$ denote unitrivalent
trees which  only differ at one  location as in Figure~\ref{fig:geomihx}. Given
an  embedded capped clasper $\Gamma_I$ of  type $I$ on a knot
$K$, then there exist capped claspers $\Gamma_H$ and $\Gamma_X$ of type $H$ and 
$X$, such that
$K_{\Gamma_I} =
\left(K_{\Gamma_H}\right)_{\Gamma_X}$.
\end{theorem} To prove this theorem, we first need the  following
\begin{proposition} \lbl{4.6} Let $K$ be an oriented knot  in a
$3$-manifold $M$, $\T$ a rooted trivalent tree, and
$E$ an edge of $\T$.
\begin{itemize}
\item[(a)] If $\G$ is a capped clasper on $K$  of type $\T$ then there is a
knot
$\tilde{K}$, and two claspers 
$\Gamma_0$ and $\Gamma_1$ of type
$\T \smallsetminus E$ on $\tilde{K}$, such  that $\Gamma_1$ is gotten from
$\Gamma_0$ by a single finger move, the guiding arc of which
corresponds to the edge $E$.
\begin{center}
\epsfig{file=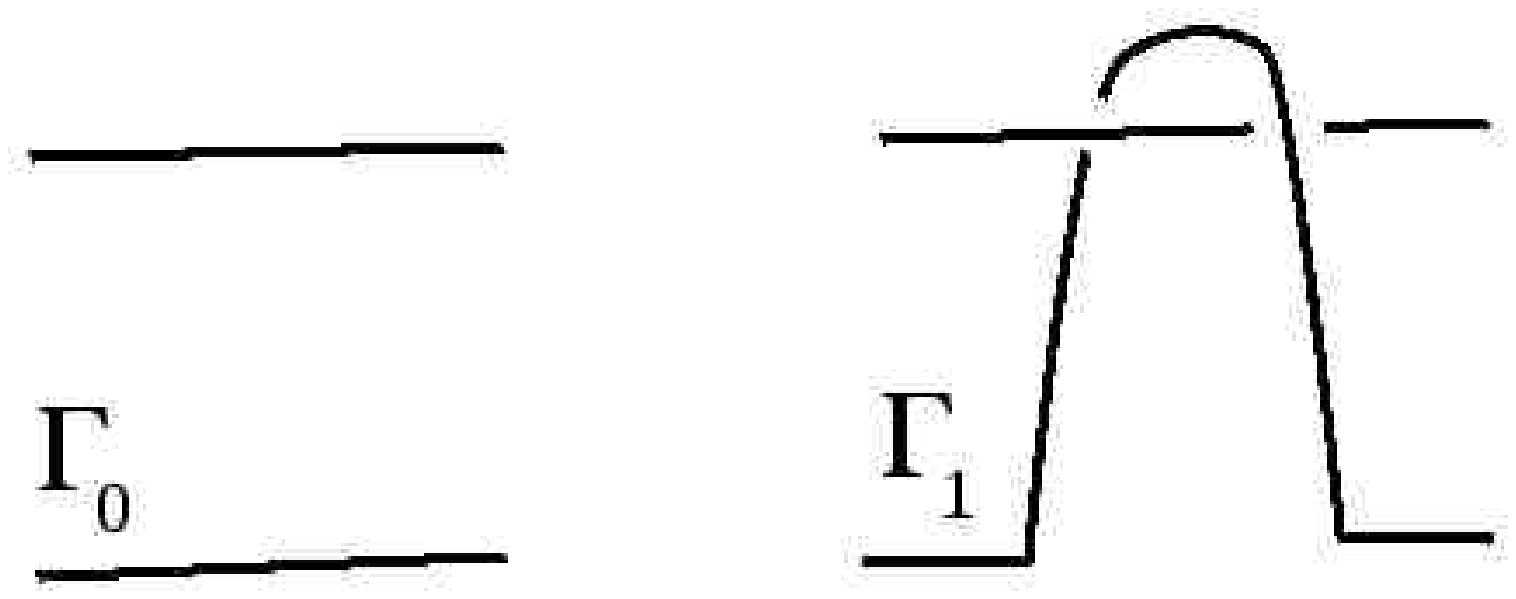,%
        height=2cm}
\end{center} and such that
$\tilde{K}_{\Gamma_0} = K$ and 
$\tilde{K}_{\Gamma_1}=K_\Gamma$.
\item[(b)] Conversely, start with  two claspers $\G_0,\G_1$ of type $\T
\smallsetminus E$ on $K$ that  differ by a finger move as above. Then there
is a clasper $\G$ of  type $\T$ such that
$$ K_{\G_1} =  (K_{\G_0})_\G
$$
\end{itemize}
\end{proposition}
\begin{proof} 

Part (a) is proven similarly to Proposition~4.6 of \cite{h2}, using a sort of 
inverse to Habiro's move 12, which is the identity in Figure~\ref{fig:newihx1}.
\begin{figure}[ht]
\begin{center}
\epsfig{file=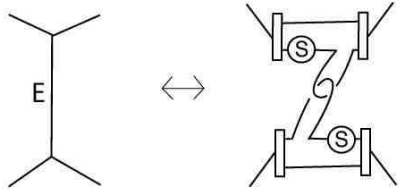,%
        height=3cm}
\end{center} 
\caption{An inverse to Habiro's move 12.}
\lbl{fig:newihx1}
\end{figure}

Now consider Figure~\ref{fig:newihx2}. One can plug either of the
two pairs of arcs (clasped respectively unclasped) on the right
of Figure~\ref{fig:newihx2} into the shaded region. After
applying Habiro's version of the  zip construction (using
claspers with boxes) as shown in the figure, one obtains a (disconnected)
clasper with boxes
$\Gamma^\prime$, and two claspers $\Gamma^\prime_0$ and
$\Gamma^\prime_1$, containing the $S$-twists. Whether one gets $\G'_0$ or
$\G'_1$ depends on what one plugs into the shaded region.
\begin{figure}[ht]
\begin{center}
\epsfig{file=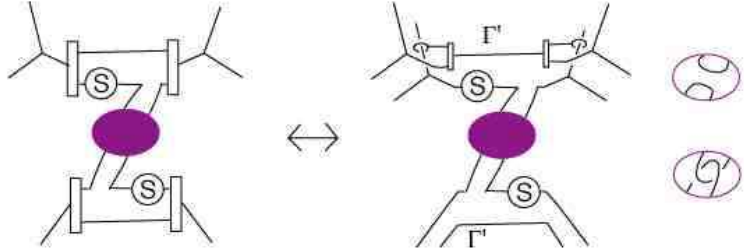, %
        height=4cm}
\end{center}
\caption{A zip move.}
\lbl{fig:newihx2}
\end{figure}

There is an important subtlety here. Surgery along a rooted
clasper (without boxes by definition) only affects the pair
$(M,K)$ inside a regular neighborhood of the clasper and
its root disk, and is fixed outside of this neighborhood.
On the other hand, for claspers with boxes, one may have to
choose many roots, modifying the pair $(M,K)$ inside a
regular neighborhood of the clasper and its root disks.
In Figure~\ref{fig:newihx2}, these added roots must include
some of the little ``lassoes'' coming out of the boxes. Hence
the clasper $\Gamma^\prime$ actually modifies
$\Gamma^\prime_1$ and $\Gamma^\prime_2$ to two claspers
$\Gamma_i= (\Gamma^\prime_i)_{\Gamma^\prime}$ for $i=1,2$.
Note that since $\G^\prime_i$ differ by a finger move, so do $\G_i$.

By the above move $K_\Gamma =
K_{\Gamma^\prime\cup\Gamma^\prime_1} =
(K_{\Gamma^\prime})_{\Gamma_1}.$  On the other hand by Habiro's
move~4, $K = K_{\Gamma^\prime\cup
\Gamma^\prime_2} =(K_{\Gamma^\prime})_{\Gamma_2}$. Thus we have
found a knot
$\tilde{K} := K_{\Gamma^\prime}$ in $S^3$ and two claspers $\Gamma_i$
which differ by a finger move
in $S^3$ and satisfy the desired identities:
$\tilde{K}_{\Gamma_0} = K$ and $\tilde{K}_{\Gamma_1} = K_\Gamma$.

Part (b) is Proposition~4.6 of \cite{h2} and  the proof is esentially the
reverse of the above  argument. \end{proof}

\begin{proof}[Proof of Theorem~\ref{Habiro}]

We  only prove the cases when the tree $I$ has at least 6 edges. The  other
case is  similar.
\begin{figure}[h]
\begin{center}
\epsfig{file=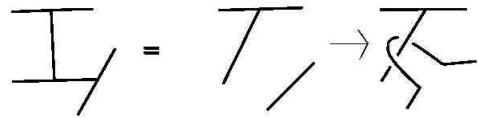,%
        height=2cm}
\end{center}
\caption{The claspers $\G_I, \G_0$ and 
$\G_1$, from left to right.  Recall
$K=\tilde{K}_{\G_0}$.}
\lbl{fig:geomihx1}
\end{figure} 

By  part (a) of Proposition~\ref{4.6} a clasper surgery on $K$
along 
$\G_I$ can be thought of as changing the clasper surgery on some knot 
$\tilde K$ from
$\G_0$ to $\G_1$ as in Figure~\ref{fig:geomihx1}. Now  apply part (b) of
Proposition~\ref{4.6} twice as  follows:
\begin{center}
\epsfig{file=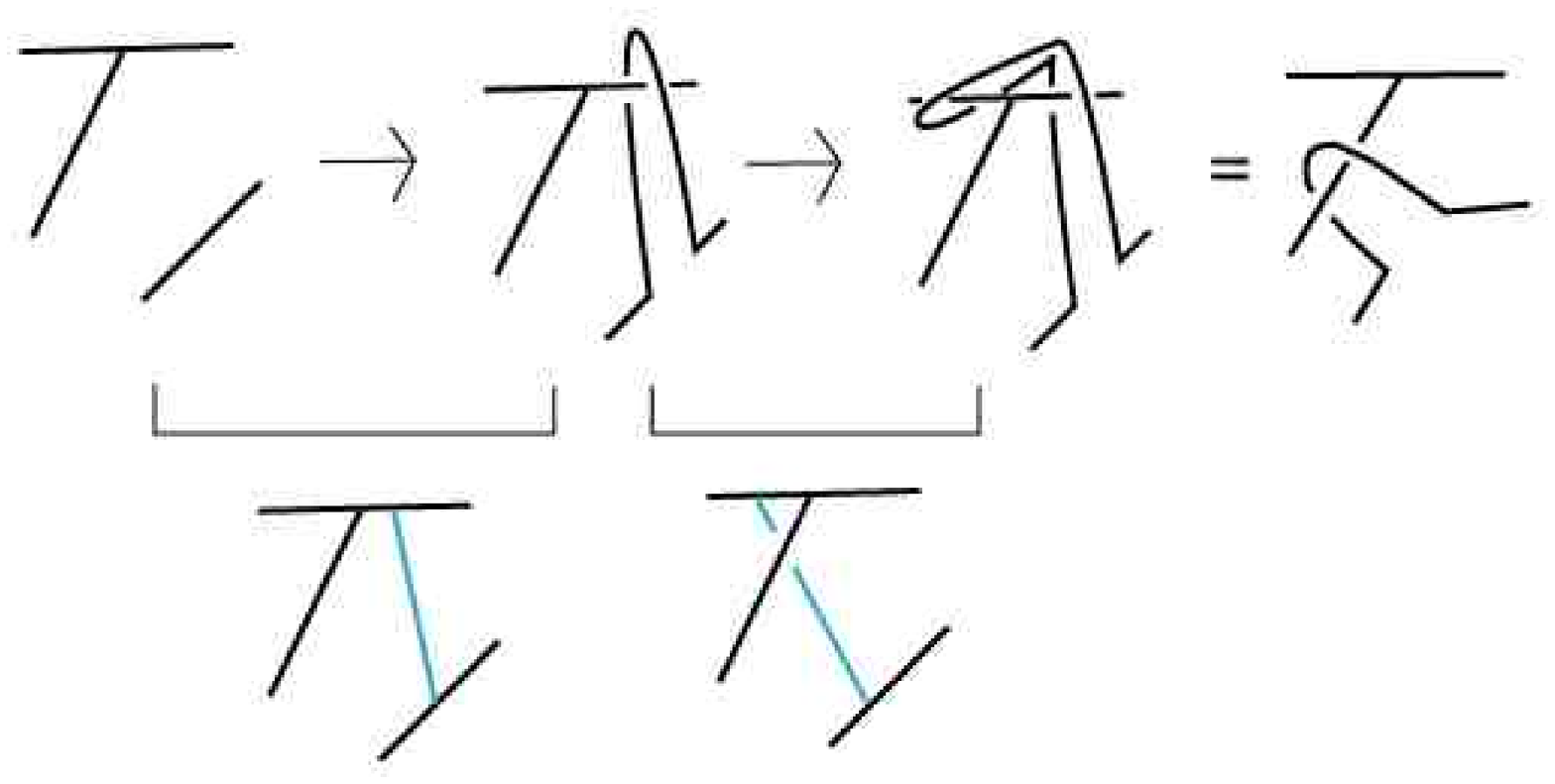,%
        height=4cm}
\end{center} This implies our  claim
$(K_{\G_H})_{\G_X}=((\tilde{K}_{\G_0})_{\G_H})_{\G_X}=\tilde{K}_{\G_1}=K_{\G_I}$.
\end{proof}

\begin{corollary}\lbl{cor:IHX}  Recall that $\h_k$ is the simplest possible
rooted tree of  class~$k$.
\begin{itemize}
\item[(a)] Capped $\h_k$-clasper surgeries  generate all capped tree clasper
surgeries of Vassiliev degree  (=class)~$k$.
\item[(b)]
$\h_k$-clasper surgeries generate all rooted  tree clasper surgeries of
Vassiliev  degree~$k$.
\end{itemize}
\end{corollary}
\begin{proof}

(a): Any  tree of class $k$ can be changed into a sequence of
$\h_k$-trees  using geometric IHX. This can be proved by introducing the
following  function on rooted class~$k$ trees $\tau$: $l(\tau )$ is the 
maximum length of a chain of edges. Given $\tau$, consider a chain of 
maximal length $c$, and suppose it misses some internal vertices. Let 
$v$ be an internal vertex of distance $1$ from $c$. Then, by  geometric IHX,
this tree can be realized as a sequence of two  trees with higher $l$:
\begin{center}
\epsfig{file=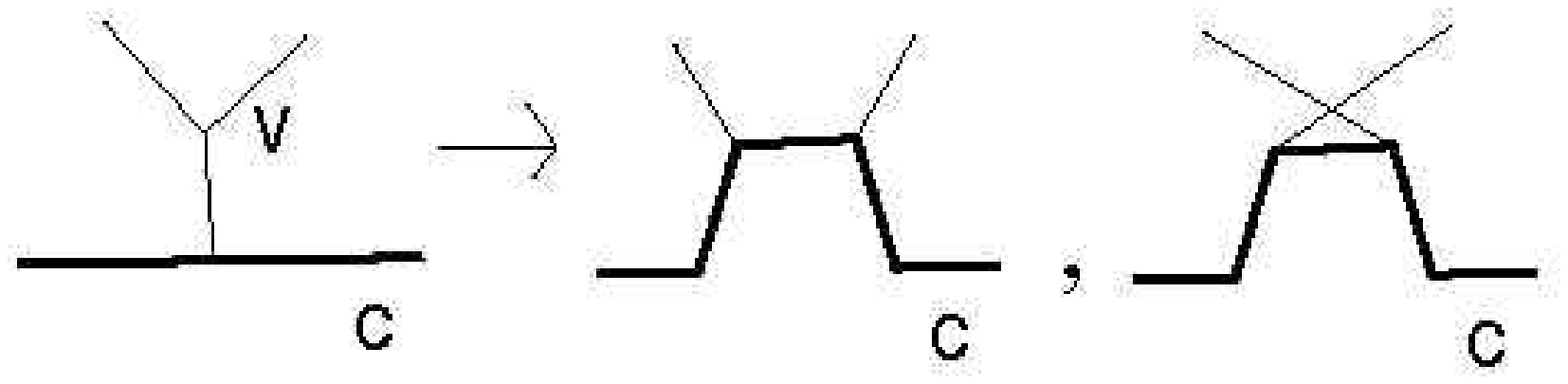,%
        height = 2.7cm}
\end{center} Hence we can keep applying IHX until we  have a sequence of
trees with maximal
$l$, which as we have seen  means that that a maximal chain hits every
internal vertex. This is  just a rooted $\h_k$-tree.

(b): Let $\C_k$ denote the set of knots  related to the unknot by
 \emph{capped} tree clasper surgeries of  Vassiliev degree $k$.  Similarly
let $\h\C_k$ denote those knots  which are related to the unknot by degree
$k$ capped tree claspers  whose tree type is that of the half grope. Define
$\rr_k$ to be those  knots related to the unknot by degree $k$ \emph{rooted}
tree clasper  surgeries,  and let $\h\rr_k$ bethe analogous object,
restricting to  half grope trees.   (By Theorems 2 and 3,
$\C_k = \F^v_k,
\rr_k = 
\F^g_k$.)

We have the following map of short exact  sequences:
\begin{equation*}
\begin{CD} 
0@>>>\h\C_k/ \h\rr_k@>>>\K/\h\C_k@>>>\K/\h\rr_k@>>>0\\ 
@VVV @VV{\text onto}V  @VV{\cong}V @VV{\text onto}V @VVV\\ 
0@>>>\C_k/\rr_k@>>> \K/ \C_k@>>> \K/\rr_k@>>>0
\end{CD}
\end{equation*}

and, by part a), the middle  map is an isomorphism. By \cite{h2}, $\K/\C_k$
is a group, a fact  which implies that everything in the above diagram is a
group (under  connected sum).
 By the $5$ lemma, the right hand map $\K/\h\rr_k \to 
\K/\rr_k$ is an isomorphism, as desired. Recall that all of the above 
quotients are defined as in the introduction, and are in particular  not just
quotient monoids. \end{proof}

The proof of  Theorem~\ref{half-gropes} is now just an application of  our
translation between gropes and claspers,  Theorem~\ref{claspers-gropes}, to
the above  Corollary~\ref{cor:IHX}.

\section{Proofs of the main results}\lbl{sec:proofs}

\subsection{Proof of Theorem \ref{claspers-gropes}}

Part (a) follows from Theorems \ref{prevthm}  and \ref{b}.

\noindent To see part (b), given a cap of a grope, this  will become a disk
bounding the corresponding leaf of the constructed  clasper, and by
definition we need to arrange that its interior is  disjoint from the
clasper. As the cap avoids the higher stages of the  grope, the only place it
might hit the clasper is along the  edge  that connects the root leaf to the
bottom stage node. Push these  intersections off the end of this edge across
the root leaf. This  introduces new (pairs of) intersections of the cap  with
the  knot, which are allowable.

Conversely, if a leaf of a clasper has a  cap, in the constructed grope the
cap will only hit the annulus part  of the bottom stage. See the discussion
after  Theorem \ref{bprime}. $\hfill\Box$

\subsection{The Zip Construction}
To prove Theorems~\ref{capped} and \ref{main} we need  a construction that will simplify a
grope cobordism to a finite sequence of moves that are simple clasper surgeries. This will be
provided in Theorem~\ref{refine} which relies on  the Habiro-Goussarov zip
construction. Habiro's version is not well suited to the present setting, since
it produces claspers with boxes, the  removal of which leads to complicated
behavior of the edges of one of the produced claspers. We state and prove a
version of the zip construction better suited to our needs. An earlier version
of this paper contained an erroneous statement of the zip contruction, which
led to an error in the statement of the original Lemma 17 which
is now replaced by Theorem~\ref{refine}. The original proof of
Theorem~\ref{capped} stays unchanged whereas the proof of
Theorem~\ref{main} now has to be supplemented by using
Corollary~4 of \cite{newconant}.

\begin{lemma}\lbl{S1} 
The following two clasper surgeries give isotopic results.
\begin{center}
\epsfig{file=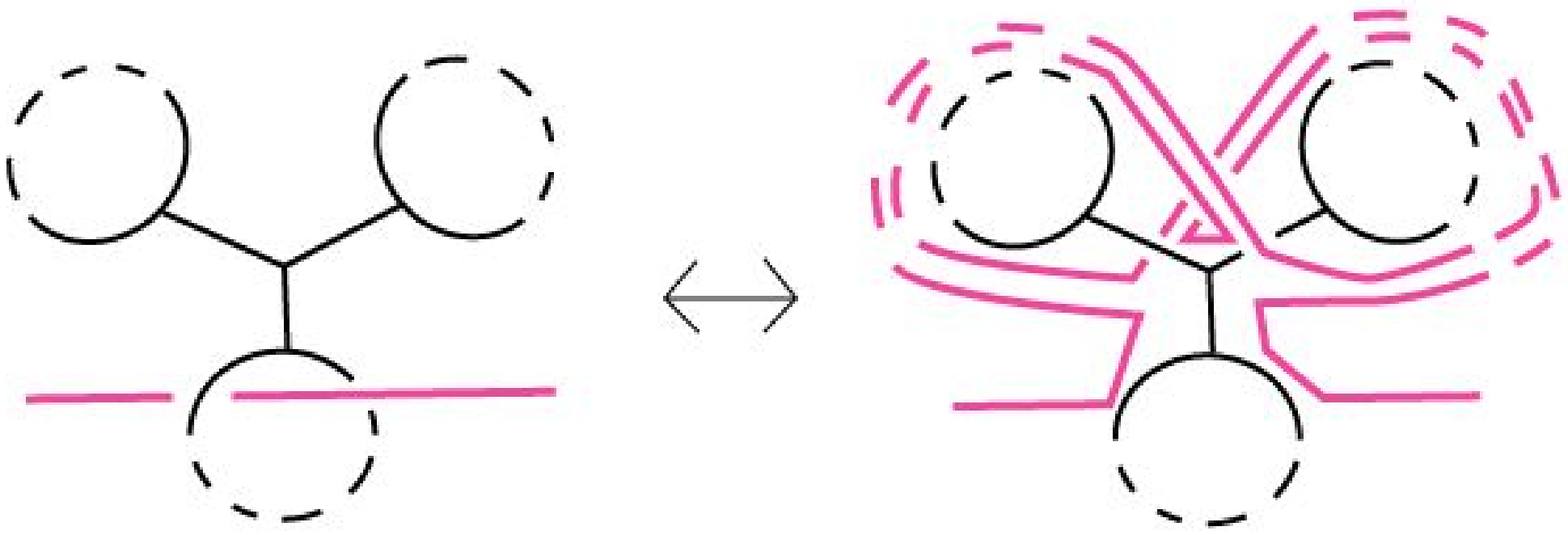,%
        height=4cm}        
\end{center}
The pictured object being slid can be part of another clasper or a strand of
the knot.
\end{lemma}
\begin{proof} 

Write out the left hand side clasper surgery as a surgery on the usual 6
component link corresponding to the Y-clasper. Then slide the visible part of the
knot or clasper over one component of the Borromean rings.
\end{proof}

\begin{corollary}\lbl{S2} Given an arc of a knot, or a piece of another
clasper that intersects a cap of a clasper $C$, then one can
slide this arc or piece of  clasper over $C$ to remove the intersection
point. That is, the slid piece lies in a regular neighborhood of $C$ minus
the leaf, and avoids any caps $C$ may have.
\end{corollary}

\begin{proof} 
Break $C$ into a union of Y-claspers and inductively
apply Lemma~\ref{S1}. 
\end{proof}

Let $L$ be a leaf of a rooted tree
clasper $C$ on a knot $K$, and let $\eta$ be a framed arc from $L$ to itself. 
Cutting the leaf along $\eta$ splits it into two halves.

\noindent {\bf Assertion:} 
Surgery on $C$ has the same effect
on $K$ as surgery on the union of  two daughter claspers $C_1$ and $C_2$,
satisfying the following properties:
\begin{enumerate}
\item $C_1$ is identical to $C$ except at $L$ where only one half of $L$ is
used.
\item The leaves of $C_2$ are parallels of the leaves of $C$ except at $L$,
where the other half of $L$ is used. The edges and nodes of $C_2$ lie in  a
regular neighborhood of $C_1$ and avoid any caps that $C_1$ may have.
\end{enumerate}
Note that in this construction one has a choice of which half of $L$ is used for
the almost-identical copy $C_1$ of $C$, and which half is used for the more
complicated daughter $C_2$.

A low degree example is shown below. 
\begin{center}
\epsfig{file=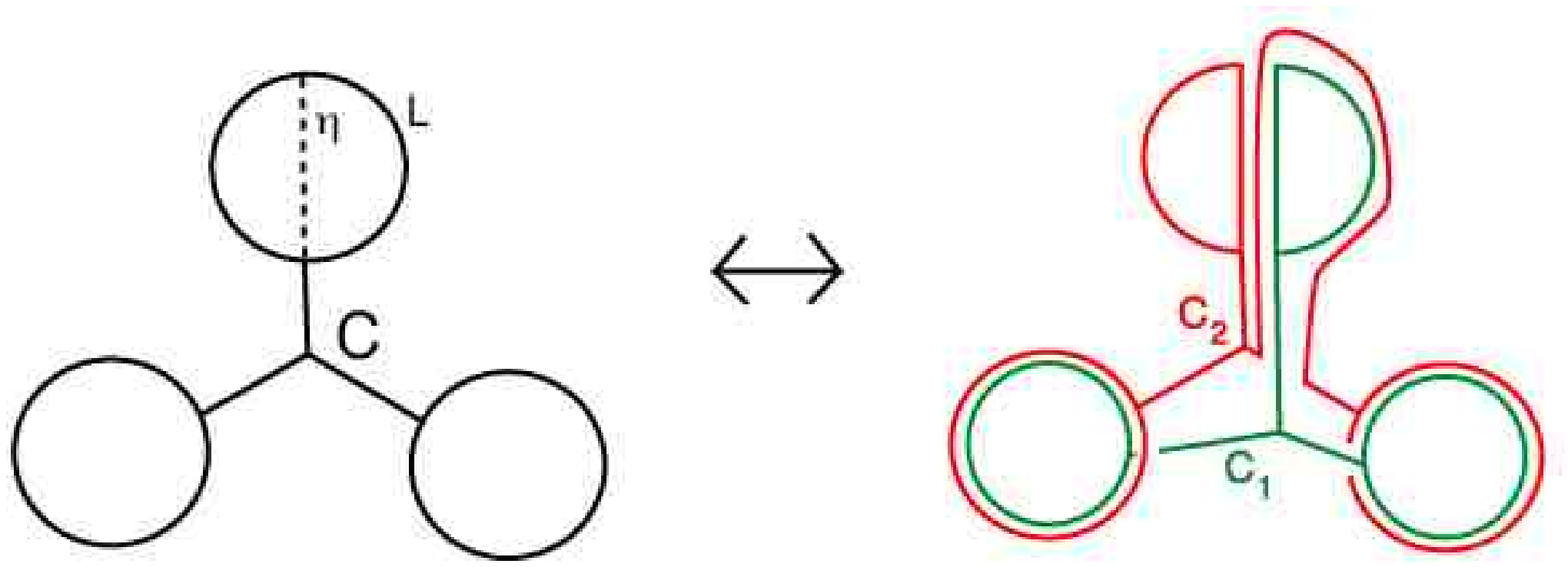, %
        height=4cm}
\end{center}
 This is in \cite{ggp}, but their
Borromean rings are oriented oppositely, so the figure should not look
identical! One can also apply the technique of Proposition 6 to obtain this
picture.

\begin{proof}[Proof of Zip construction, i.e. of the assertion above]
The statement follows from the following more general statement: 
Inside a regular neighborhood, $N$, of $C\cup\eta$, there are two claspers
$C_1$ and $C_2$ as above, such that $N_C$ is diffeomorphic
rel boundary to $N_{C_1\cup C_2}$. Notice that since $C_2$ avoids any caps 
that $C_1$ may have, it in particular avoids the root leaf. 

We proceed by induction, the picture above serving as the base case. In
the pictures that follow, the thicker lines denote a regular
neighborhood of a clasper. To induct, we break the clasper $C$
into a union of two simpler claspers as follows:
\begin{center}
\epsfig{file=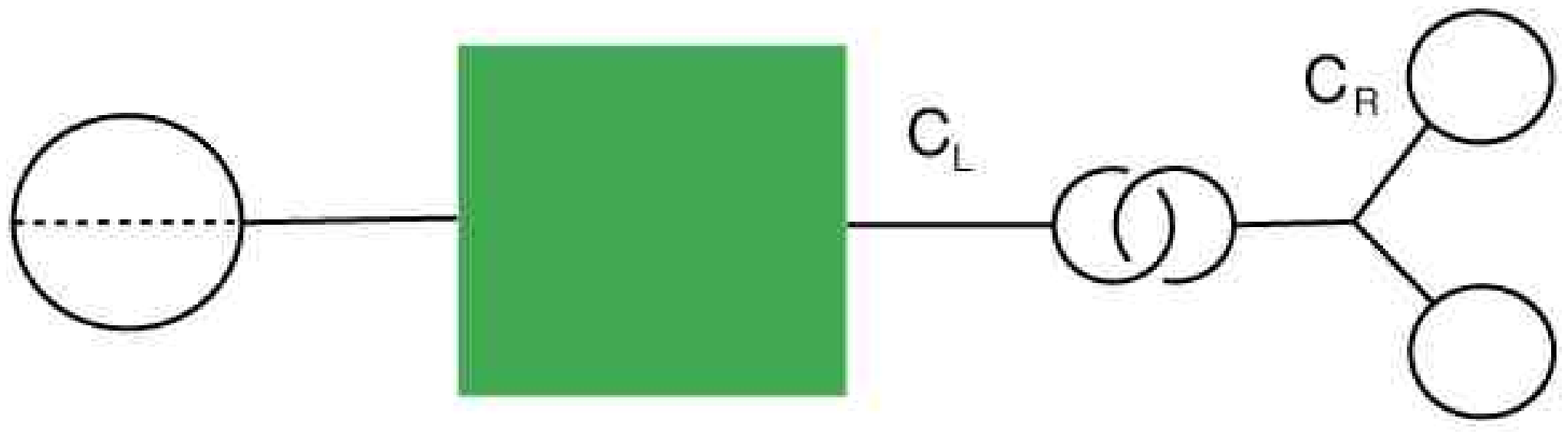, %
        height=3cm}
\end{center}
The big box is a pictorial convenience to represent an arbitrary clasper.
Inductively we get the following picture:
\begin{center}
\epsfig{file=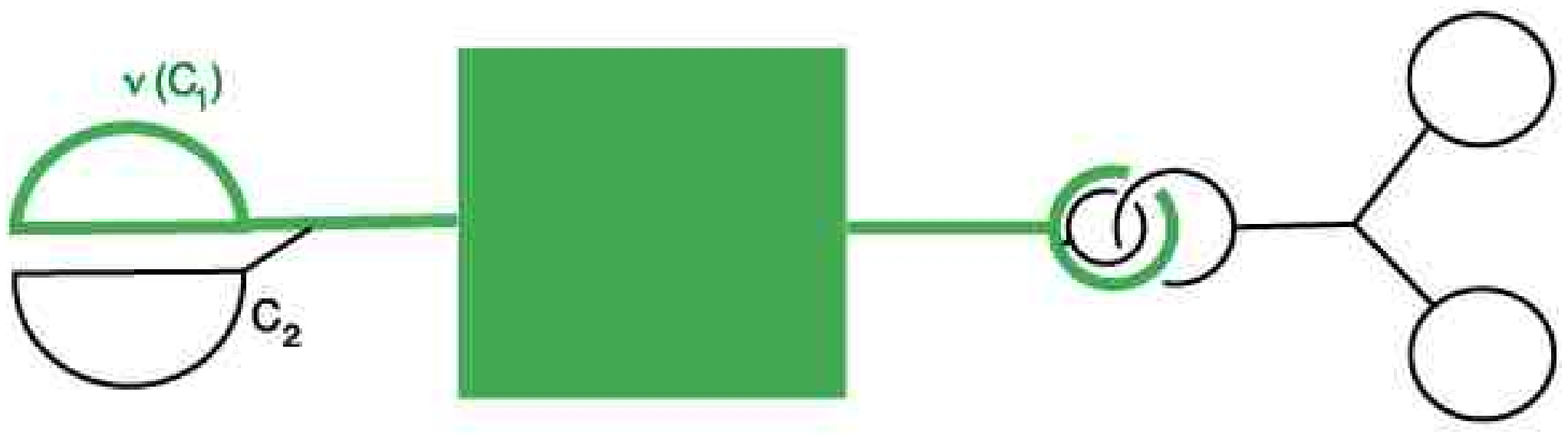, %
        height=3cm}
\end{center}
Then using the base case on the left leaf of  the right-hand
clasper, we obtain
\begin{center}
\epsfig{file=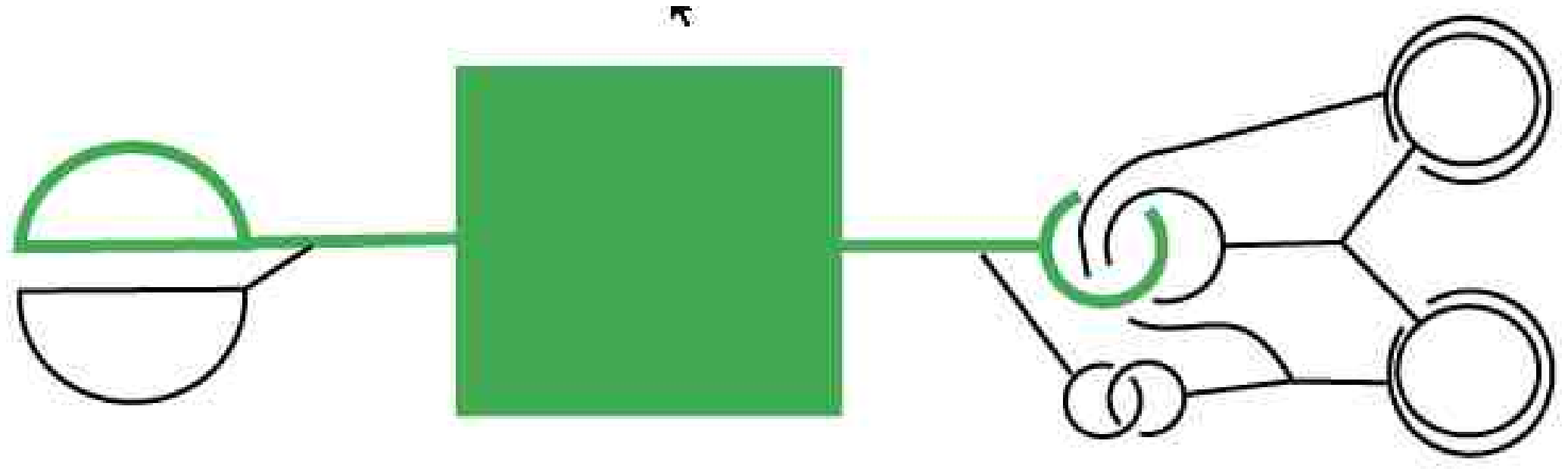, %
        height=3cm}
\end{center}
By Corollary~\ref{S2} applied to the grey leaf on the right and
a cancellation of the bottom Hopf pair, we get
\begin{center}
\epsfig{file=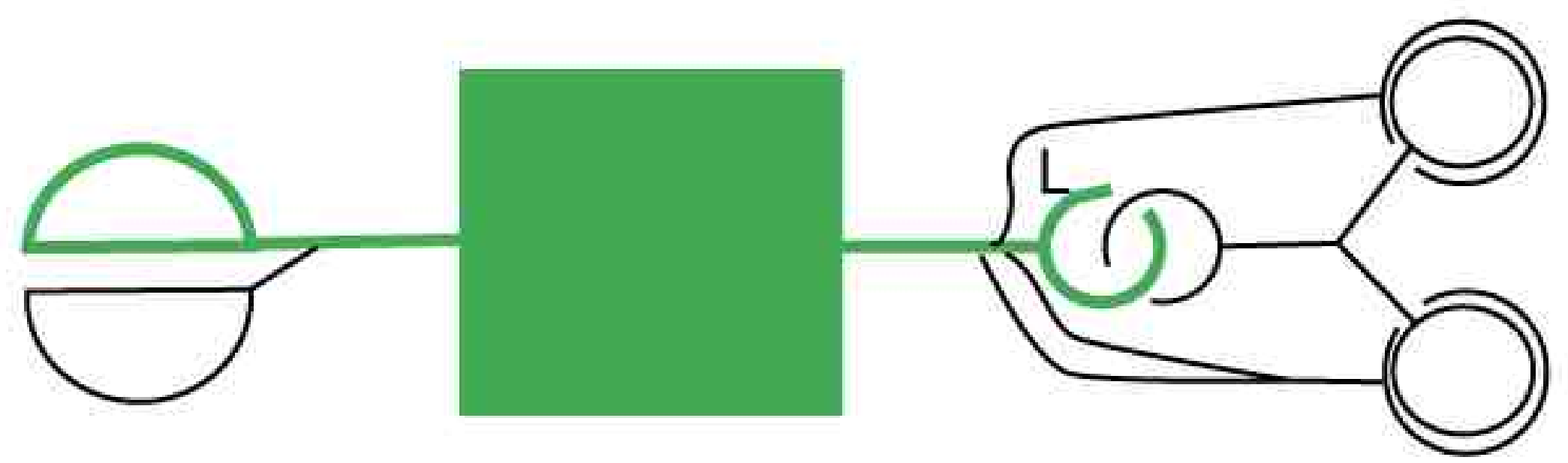, %
        height=3cm}
\end{center}
Next we would like to cancel the grey-black Hopf pair above. This requires some
care because parts of $C_2$ run parallel to the grey leaf $L$.
However, in our construction, $C_2$ avoids the caps of $C_1$.
Thus we can split the regular neighborhood of $L$ apart into
the leaf, plus a parallel copy of that leaf through which other
claspers wander:
\begin{center}
\epsfig{file=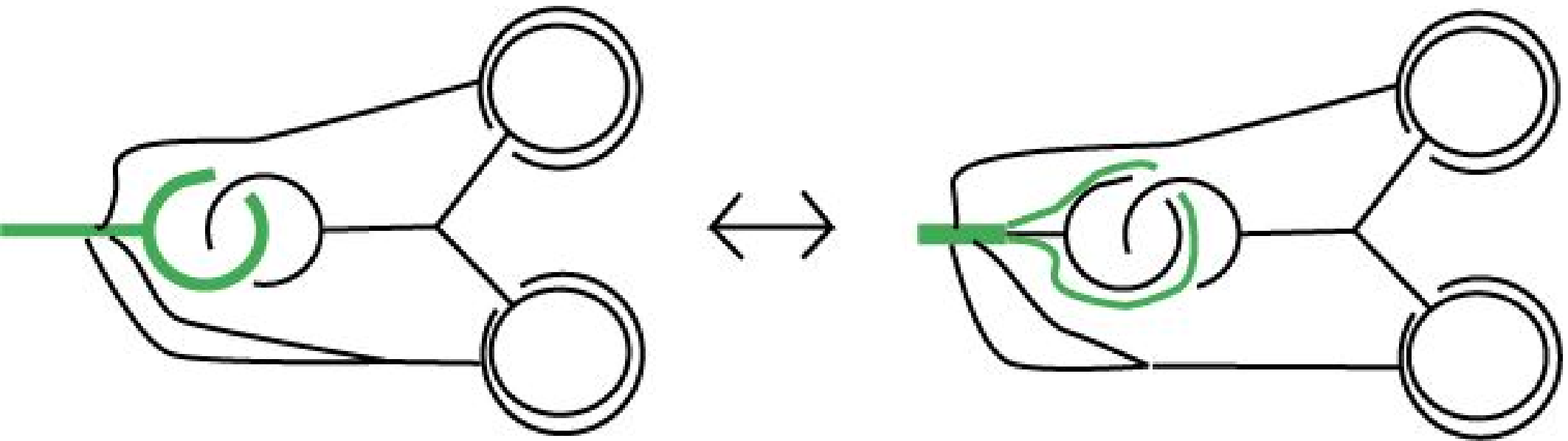, %
        height=3cm}
\end{center}
After that we apply a sequence of Corollary~\ref{S2} moves to obtain a clean
Hopf pair that can be cancelled. In the figure below we also push some black
arcs into the grey area which after all only represents some neighborhood
of the clasper:
\begin{center}
\epsfig{file=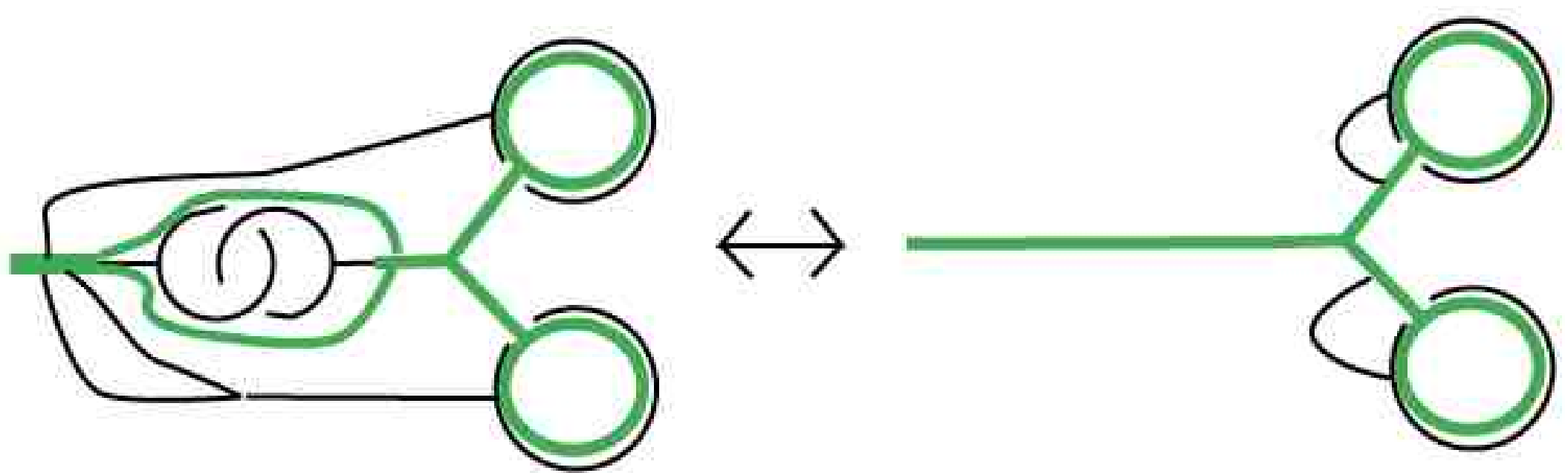, %
        height=3cm}
\end{center}
Thus we have finished the inductive step.
\end{proof}

\subsection{Simplifying a grope cobordism}

\begin{lemma}\lbl{R2} Let $C$ be a rooted tree clasper of type $\T$ 
 with a leaf $L$ bounding a disk that only intersects edges of $C$ (and is
disjoint from the knot $K$). Then the surgery on
$C$ may be realized as a sequence of clasper surgeries along claspers
$C_1,\dots,C_n$ which come in two types:
\begin{itemize}
\item[(a)] $C_1$ is identical to $C$, except that the leaf $L$ is replaced by
a leaf that has a cap. In particular, $C_1$ has type $\T$.
\item[(b)] $C_i$, for $ i>1$, have type $\T^\prime$, where $\T^\prime$ is the
tree formed from
$\T$ by gluing a ``Y'' onto the univalent vertex representing $L$. In
particular, the degree of $\T^\prime$ is bigger than that of $\T$.
\end{itemize}
\end{lemma}

\begin{proof}
Push each intersection point of an edge with the given disk bounding $L$  out
toward the other leaves, using little fingers following the spine of the clasper
$C$.  Each such finger splits into two at a trivalent vertex of $C$, and stops
right before a leaf (which is necessarily distinct from $L$). This describes a
new disk $D$ bounding $L$ which has the property that
on each edge $E_i$ incident
to a leaf $L_i\neq L$ there are several parallel sheets of $D$ being punctured by
$E_i$ (and there are no intersections of $D$ with edges other than $E_i$). If the
leaf
$L_i$ happens to be the root leaf, we push these sheets over the cap of $L_i$,
introducing intersections with the knot, but eliminating the intersections with
$E_i$. If
$L_i$ is not the root, we add a series of nested tubes that go around $L_i$,
trading the intersections with $E_i$ for genus on $D$.

Thus $L$ now bounds an embedded surface which intersects $K$ but is disjoint from
the clasper $C$. We perform the zip construction on $L$ to segregate the
knot intersections, where the first daughter $C_1$ will inherit the half of $L$
bounding a disk intersecting the knot. This first daughter is of type (a). The
second daughter has the leaf coming from the half of $L$ bounding a surface
disjoint from the clasper $C_2$. Converting  the clasper to a grope we get a
grope of tree type
$\T$ whose tip corresponding to $L$ bounds a surface disjoint
from the grope. Hence we really have a grope of increased
class, but it has high genus at the tip $L$. Proposition 6 now
yields a sequence of cobordisms of type
$\T^\prime$ as claimed.   
\end{proof}

The  following cleaning up procedure is the heart of this section. It is in
spirit similar to the procedure described in section 4.3 of
\cite{ggp}. There the authors work in the context of Goussarov's finite type
theory (using alternating sums to define a filtration on the span of all knots).
Here we need to strictly work with clasper moves on knots, there are no linear
combinations that can help with cancellations. Therefore, the geometric arguments have to be
much more subtle.

\begin{theorem} \lbl{refine} Let $\T$ be a rooted  trivalent tree. We can
realize any $\T$-grope cobordism in $S^3$ by a  sequence of clasper surgeries
each of which either has higher grope degree than the  original, or is a
$\T$-clasper surgery  which has tips of the following  form:
\end{theorem}

\begin{center}
\epsfig{file=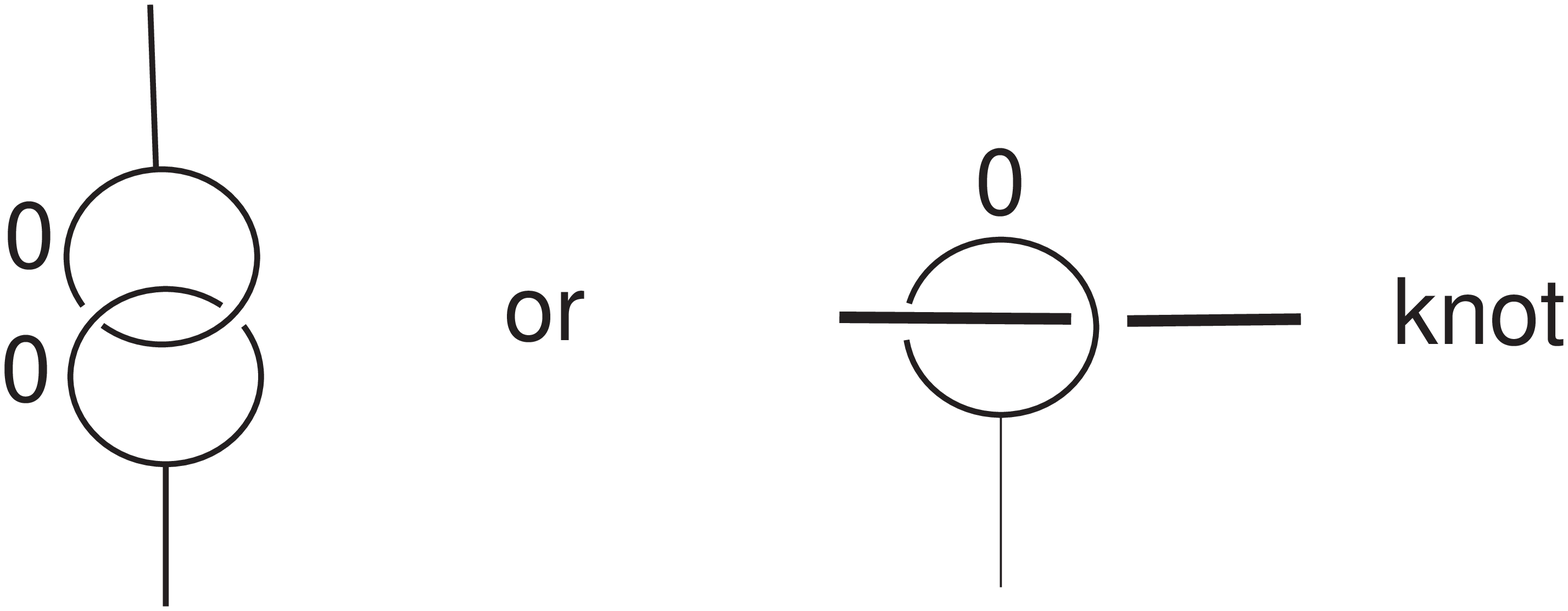,%
        height = 3cm}
\end{center}

\begin{proof}

By proposition 
\ref{reducegenus} we may assume that all surface stages of the given grope are
of genus one.  Such a grope cobordism corresponds to a  $\T$-clasper surgery,
which we proceed  to simplify.

\vspace{3mm}
\noindent {\bf Step 1}: First we make the leaves $0$-framed. This is 
accomplished using the following simple observation. Suppose $x$ and 
$y$ represent symplectic basis elements on a punctured genus one  surface
embedded in $S^3$. These have framings $\sigma(x),\sigma(y)$,  the diagonal
terms of the Seifert matrix. There is also the  intersection pairing
$\I:H_1(F)\otimes H_1(F)\to \Z$. By assumption 
$\I(x,y)=1$. The  formula
$$\sigma(a+b)=\sigma(a)+\sigma(b)+\I(a,b)
$$ implies that if 
$\sigma(x) = n$ and $\sigma(y) = 0$, then
$\sigma (x-ny) = 0$. By  Dehn twisting one can represent $x-ny, y$ by embedded
curves meeting  at a point. So $x-ny,y$ represent a $0$-framed basis of $F$.
In  particular, suppose $F$ is a surface stage of the grope for which $x$  is a
tip, and $y$ bounds a higher surface stage. Then $\sigma(y) = 0$  and we can
let $x-ny$ be the tip in place of $x$. This takes care of  all possibilities
except the case when $x$ and $y$ are both tips of  the grope which have
nonzero framings. Here we perform some  sleight-of-hand using claspers. 
Convert the grope to a clasper
$C$.  Then insert a Hopf-linked pair of leaves on the edge incident to $y$. 
This disconnects the clasper into two pieces $C_x,C_y$ as  in
Figure~\ref{fig:zero}.
\begin{figure}[ht]
\begin{center}
\epsfig{file=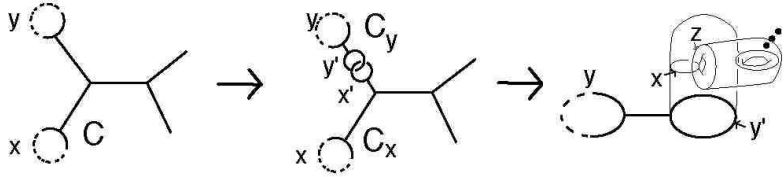,%
        height = 2.6cm}
\end{center}
\caption{Making $x$ 
$0$-framed.}
\lbl{fig:zero}
\end{figure}

The tips $x$ and $y$ each  lie on exactly one of these claspers. The other
leaf $y'$ of $C_y$  bounds a grope
$\tilde{G}$ gotten from $C_x$, by considering $x'$ as  the root leaf.
$z$ is the curve on the bottom stage of $\tilde{G}$  which bounds the next
surface stage, as pictured. By changing $x$ to 
$x-nz$ as before, we convert the tip $x$ of 
$\tilde{G}$ to  a zero-framed tip. Changing $\tilde{G}$ to a clasper
$C^\prime_x$ by 
our procedure, we again have the clasper
$C_y$ with the leaf $y'$  Hopf-linking the root $x'$ of $C^\prime_x$. 
Convert this back to an  edge to achieve a clasper of the same type as $C$,
but with one  more tip zero framed. This clasper may be converted back to a
grope  if we  wish. Notice that under our grope-clasper correspondence, the 
framings of tips  (leaves) do not change.  Do this until all tips are 
zero-framed.

\vspace{3mm}
\noindent {\bf Step 2}: Next we make the leaves  unknotted. It is an exercise to
prove that there is a set of arcs  from a knotted leaf to itself, such that
cutting along these arcs  yields a collection of unknots. Hence, given a knotted
leaf, one can  apply the zip construction to such a set of arcs, thereby
reducing  the number of knotted leaves in each resultant clasper. Repeat this 
procedure until you have a set of claspers with unknotted  leaves.

Note that we have now proved that any 
$\T$-clasper surgery can be reduced to a sequence of $\T$-clasper  surgeries,
each of which has $0$-framed leaves bounding disks. To continue, we need
to clean up the intersection pattern of the disks.  By pushing fingers of disks
out to the boundary, one may assume each  pair of disks intersects in clasp
singularities; i.e. the  intersection pattern on each disk is a set of arcs
from interior  intersections with the clasper to the boundary of the disk. 
Secondly, we eliminate triple points. After we did the first step,  there is
a triple point which is connected by a double point arc to  the boundary of
one of the disks, such that there are no intevening  triple points. Push a
finger of the disk which is transverse to this  arc along the arc and across
the boundary.  Repeat this until all  triple points have been removed. This
homotops the disks into a  position such that the intersection pattern
consists of disjoint  clasp singularities.

\vspace{3mm}
\noindent {\bf Step 3}:
We now start with a clasper $C$ which has $0$-framed leaves bounding disks
$D_i$ with only clasp intersections between each other.  
In addition, the disks $D_i$ may have several types of intersections with $C$
and the knot $K$, which we proceed to organize. Note that our theorem states
that, modulo higher grope degree, we can reduce to only two types of
singularities for the
$D_i$: Either there is a single clasp (and no other intersections with $C$ or
$K$), or there is a single intersection with the knot
$K$ (and no intersections with $C$). We call such disks {\em good} for the
purpose of this proof. The bad disks fall into several cases which we will
distinguish by adding an index to the disk
$D$ which explains the failure from being good. The cases are as follows, where we
list {\em exactly} the singularities of the disk, so unmentioned problems do not
occur. 

\noindent If a disk $D$ has\dots 

\noindent \dots intersections with edges of $C$, we call it $D_E$.

\noindent\dots  more then one intersection with $K$, we call it $D_K$.

\noindent\dots  has more than one clasp, we call it $D_{Cl}$.

\noindent\dots  intersections with edges of $C$ and with $K$, we call it
$D_{E,K}$.
 
\noindent \dots  intersections with edges of $C$ or with $K$, and has clasps, we
call it $D_{EK,Cl}$.
\vspace{3mm}

Just to be clear, the cases  $D_E$, $D_K$ and $D_{E,K}$ above
represent disks without clasps, whereas $D_{Cl}$ has no intersections with edges of
$C$ or with $K$.

It  is clear that these cases represent all possibilities for a bad disk. Recall
that a disk $D$ was called a cap if it is embedded disjointly from $C$. In our
notation, this means that a cap is either bad of type $D_K$ (more than one
intersection with $K$), or it is good (exactly one intersection with $K$). We
ignore the case of a cap without intersections with $K$ since then the surgery
on the clasper has no effect on $K$.

We now introduce a complexity function on claspers with given disks $D_i$ as
above.  It is defined as a quintuplet $(c_1,c_2,c_3,c_4,c_5)$ of
integers $c_i$, ordered lexicographically. The $c_i$ are defined as
follows:
\begin{itemize}
\item $c_1$ is minus the number of disks $D_i$ which are caps. 
\item $c_2$ is the total number of intersections of the knot with caps $D_i$.
\item $c_3$ is the total number of clasps.
\item  $c_4$ is the number of bad disks of type $D_{E,K}$.
\item  $c_5$ is the number of bad disks of types $D_{EK,Cl}$.
\end{itemize}

The proof proceeds by using the zip construction to split a bad disk of a clasper
into two daughters. In each of the five cases given below we
check that both daughter claspers have either smaller
complexity or higher grope degree, so they are ``cleaned up''. 
The five cases can be applied in an arbitrary order and they
are performed as long as there is a bad disk on a daughter clasper (where we
do not work on claspers of higher grope degree). Since each $c_i$ is bounded
below, this cleaning up process must terminate. This can only happen if all
disks are good (or the clasper has higher grope degree), which is the statement
of our theorem.  

We now describe the five cases of the cleaning up process. In
each case the label says which bad disk is being split, then we
have to specify the splitting arc and the order of the daughter
claspers. 

\noindent (E) \\
Suppose there is a bad disk of type $D_E$. By Lemma~\ref{R2},
this splits into a daughter clasper $C_1$ of the same degree but with an extra
cap, and into a sequence of claspers of higher grope degree.  For $C_1$ the
number $c_1$ is reduced.

\noindent (K)\\
 Suppose there is bad disk of type $D_K$.
Split along an arc that divides the intersections  with $K$ into two smaller
sets. Each daughter clasper inherits a cap with fewer intersections, so $c_2$
goes down for both daughters (whereas $c_1$ is unchanged).

\noindent (Cl) \\
Suppose there is bad disk of type $D_{Cl}$. Draw an arc along
the disk  separating the clasps into two smaller groups. The zip construction
produces two daughter claspers $C_1$ and $C_2$ for which $(c_1,c_2)$ are
preserved. To calculate the change in $c_3$ we need only consider the leaves of
$C_1$ and
$C_2$ as $c_3$ does not see knot or edge intersections. The leaves of $C_i$
differ from those of $C$, only by cutting off part of the leaf we are
splitting along. By construction, this has fewer clasps, i.e. $c_3$ is reduced
for both daughters $C_i$.

\noindent (E,K) \\
Suppose there is bad disk of type $D_{E,K}$. Split along an arc
separating the two types of intersections, such that $C_1$ inherits the part of
the leaf with just edge intersections. Since the intersection pattern for $C_1$
is just a subpattern of the original, the entire complexity function cannot
increase. But
$c_4$ clearly decreases for $C_1$ because a new disk with only edge intersections
has been created. On the other hand $C_2$ has a new cap, so $c_1$ decreases for
it.

\noindent (EK,Cl) \\
Suppose there is bad disk of type 
$D_{EK,Cl}$.  Split along an arc which separates the clasps from the other types
of  intersections. Split in such a way that $C_1$ inherits the part of the leaf
which has the clasps. Now $(c_1,c_2)$ is preserved in $C_1$.  The cut leaf now
has only clasp intersections, and since the  intersections of the disks of
$C_1$ with everything are decreased,  new disks with both clasp and other
types of intersections are not created. Hence $c_5$ decreases for $C_1$. Now we
analyze $C_2$.  Since $(c_1,c_2)$ can only go down when we split, it suffices to
show that $c_3$ decreases. This follows by the same argument as case (Cl). 

We note that in the above five cases, when we split along a
disk, the caps away from the split disk are preserved, as are
the number of intersections of the knot with these caps.
Furthermore, in the first daughter clasper $C_1$ the four
complexity functions
$c_1, c_3,c_4,c_5$ must each stay the same or go down, because
the intersection pattern of $C_1$ is just a subpattern of the
one for the original clasper. The number $c_2$ can only
increase during an $(E)$-move, but then $c_1$ goes down for the
first daughter $C_1$ (and $C_2$ has higher grope degree).

The intersection pattern for $C_2$ changes in a more complicated way. The  first
problem is that it sits on a different knot: the knot modified by $C_1$,
which adds intersections of the knot with the disks $D_i$.  (We are applying
$C_1$ and $C_2$ sequentially!) The second problem is that the edges of $C_2$
wander around inside a neighborhood of $C_1$ and add intersections as well.
Therefore, the complexities $c_4$ and $c_5$ may increase from $C$ to $C_2$ in
all moves above, except for $(E)$.

 We summarize the information of these moves in
the following table.  Observe that performing
a move always implies a reduction of the relevant complexity $c_i$, which
we have written first in its row. Other complexities may or may not
increase, and in some cases they actually decrease. In that sense the table
contains the worst case scenario for the complexities $c_i$ of the two daughter
claspers. The notation $c_i\up$ means that $c_i$ may increase (which is
bad), whereas $c_i\down$ is the good case where the complexity definitely
decreases. Unmentioned complexities $c_i$ either stay unchanged or decrease.

\vspace{3mm}
\begin{tabular}{|c||l|r|}
\hline
Move & First daughter & Second daughter \\
\hline\hline
(E) & $c_1\down, c_2\up$ & higher grope degree  \\
\hline
(K) & $c_2\down$  & $c_2\down, c_4\up, c_5\up $ \\
\hline
(Cl) & $c_3\down$  & $ c_3\down ,c_4\up, c_5\up$ \\
\hline
(E,K) & $c_4\down $ & $ c_1\down,c_2\up,c_4\up, c_5\up  $ \\
\hline
(EK,Cl) & $c_5\down$ & $c_3\down ,c_4\up, c_5\up $ \\
\hline
\end{tabular}

\vspace{3mm}
\noindent We see from this worst scenario table that for all the
five moves the total complexity goes down for both daughter
claspers (or the grope degree increases). This completes our
argument.
 \end{proof}

\subsection{Proof of  Theorem~\ref{main}}
Consider a grope cobordism of tree type $\T$  (and
class $c$) between two knots
$K_1$ and $K_2$ in 3-space. The preceding Theorem~\ref{refine} allows us  to
reduce each of these  to a sequence of
$\T$-clasper surgeries with leaves of only two possible good types,
together with claspers of higher degree. Applying Theorem~\ref{refine} again to
these higher degree terms, and iterating, we obtain a sequence of claspers of
degrees $c$ to
$2c$ each of which has only the two good types of leaves, together with some
claspers of degree~$(2c+1)$. By Theorem~3 of
\cite{newconant} a rooted clasper $C$ of degree~$(2c+1)$  preserves
Vassiliev-Goussarov equivalence of degree~$c$. Then, by the main
theorem of
\cite{h2}, surgery on $C$ can be realized as a sequence of simple tree clasper
surgeries of degree $c$.
Recall that a {\em simple} tree clasper in Habiro's sense has by definition only
the simplest type of leaf, namely bounding a cap which intersects the knot once.
This is one of the good leaf types from Theorem~\ref{refine}.

Thus we get a sequence
of tree claspers in degrees~$c$ to $2c$ each of which only has the two good types
of leaves. For each such tree clasper, convert the Hopf-linked pairs  of leaves
to edges (or half-twisted edges). Observe that the  resulting graph claspers are
simple, i.e. they are capped and the  knot intersects each cap in exactly one
point. Let
$G$ be the graph type of one of these  simple claspers. Then the loop degree
$\ell(G)$ is the number of  Hopf-linked pairs of leaves because we started
with a tree $T$ and  glued up pairs of tips. Each such gluing reduces the
number of  vertices by two and hence the grope degree is unchanged from $T$
to $G$:
$$ 
g(G)=\ell(G)+v(G)=g(T)=v(T)\in[c,2c].
$$ 
This implies that 
$[K_1]=[K_2]\in \K/\F^g_c$ because by definition the equivalence  relation
corresponding to $\F^g_c$ is generated by simple clasper  surgeries of grope
degree~$\geq c$.

Conversely, if $C$ is a simple  clasper of type $G$ (and grope degree~$c$),
then we can convert 
$\ell(G)$ edges into Hopf-linked leaves as in Figure~\ref{fig:morse1}  to
obtain a simple tree clasper of the same grope degree, which now has class $c$.
Picking any leaf as the root, our main  construction,
Theorem~\ref{claspers-gropes}, gives a grope cobordism  of class $c$.
$\hfill\Box$

\subsection{Proof of  Theorem~\ref{capped'}}
By Theorem~\ref{claspers-gropes}, two knot types are capped grope cobordant of class~$c$ if
and only if they are related by a sequence of capped tree clasper surgeries of class (or 
Vassiliev degree)~$c$. Applying the  algorithm of Theorem~\ref{refine} (case (K) is all that
is needed) to a capped tree clasper, we get a  sequence of {\em simple} tree claspers of the
same type (and hence class). This uses the fact  that the algorithm never introduces
intersections between a cap and the clasper. Thus two knots which are capped grope cobordant
of class~$c$ do represent the same element in $\K/\F^v_c$ (the equivalence relation
generated by simple clasper surgeries).

Conversely, if two knots represent the same element in $\K/\F^v_c$, then
by Habiro's main theorem  they are also related by a sequence of simple {\em tree}
clasper surgeries of class~$c$, and thus they are capped grope cobordant
of class~$c$.
$\hfill\Box$

\subsection{Proof of Theorem \ref{concordance}}\lbl{sec:4D}

Turn the simple clasper $C$ into a  tree clasper by converting some edges
into Hopf-linked pairs of  leaves. Notice that all the resulting leaves bound
disks into the  complement of $L$. Picking a root of
$C$, and hence the corresponding  component $L^0$ of $L$, this gives a
$3$-dimensional grope cobordism  between $L$ and $L_C$. Since $\ell(C)\geq 1$
there is one tip which  bounds a cap into the complement of
$L$. Push the interior of this  cap slightly up into $S^3\times I$. Now extend
$L$ by annuli up to 
$\R^3\times{1}$. These annuli miss the pushed-up cap by construction.  The
result is an embedded grope connecting $L^0_C$ and $L^0$  in
$\R^3\times [0,1]$,
   with one tip bounding an embedded cap. The  usual procedure of iterated
surgery on this cap produces an annulus  which is disjoint from the straight
annuli connecting the other  component of $L_C$ and $L$. Thus we have
constructed a concordance,  which at closer inspection turns out to be a
ribbon concordance. This  follows from the fact that the only nontrivial
parts come from copies  of the cap which was pushed up from $\R^3$ into $\R^3
\times [0,1]$.  Hence reading from $L_C$ to $L$, the concordance has only
local  minima and saddles, but no local  maxima. 
$\hfill\Box$

\end{document}